\documentclass[a4paper,12pt]{article}
\usepackage{amssymb,amscd,amsfonts,amsbsy}
\usepackage{enumerate}
\usepackage{epsfig}
\input euscript.sty
\usepackage{color}

\def\ring{\mathaccent"0017 }

\newcommand{\RR}{{\mathbb{R}}}
\newcommand{\NN}{{\mathbb{N}}}
\newcommand{\ZZ}{{\mathbb{Z}}}

\newcommand{\PP}{{\mathbb{P}}}

\newcommand{\eps}{\varepsilon}
\newcommand{\po}{\partial\Omega}
\def\ring{\mathaccent"0017 }

\newcommand{\bp}{\noindent {\it Proof}.\,\,}
\newcommand{\ep}{\hfill$\Box$ \vskip 0.08in}

\newtheorem{proposition}{Proposition}[section]
\newtheorem{theorem}[proposition]{Theorem}
\newtheorem{lemma}[proposition]{Lemma}
\newtheorem{corollary}[proposition]{Corollary}

\begin{document}

\title{Boundedness of the gradient of a solution and \\Wiener test of order one for the
biharmonic equation}


\author{Svitlana Mayboroda (svitlana@math.brown.edu),\\
 Vladimir Maz'ya (vlmaz\@@math.ohio-state.edu)\thanks{
2000 {\it Math Subject Classification:} 35J40, 35J30, 35B65.
\newline {\it Key words}: Biharmonic equation, Dirichlet problem,
general domains.
\newline The second author is partially supported by NSF grant DMS 0500029.}}

\date{ }

\maketitle

\begin{abstract} The behavior of solutions to the biharmonic
equation is well-understood in smooth domains. In the past two
decades substantial progress has also been made for the polyhedral
domains and domains with Lipschitz boundaries. However, very
little is known about higher order elliptic equations in the
general setting.

In this paper we introduce new integral identities that allow to
investigate the solutions to the biharmonic equation in an {\it
arbitrary} domain. We establish:\\
(1) boundedness of the gradient of a solution in any three-dimensional domain;\\
(2) pointwise estimates on the derivatives of the biharmonic Green
function;\\
(3) Wiener-type necessary and sufficient conditions for continuity
of the gradient of a solution.
\end{abstract}

\section{Introduction}
\setcounter{equation}{0}

The maximum principle for harmonic functions is one of the
fundamental results in the theory of elliptic equations. It holds
in arbitrary domains and guarantees that every solution to the
Dirichlet problem for the Laplace equation, with bounded data, is
bounded. In 1960 the maximum principle has been extended to higher
order elliptic equations on smooth domains (\cite{Agmon}), and
later, in the beginning of 90's, to three-dimensional domains
diffeomorphic to a polyhedron (\cite{KMR}, \cite{MR1}) or having a
Lipschitz boundary (\cite{PVmax}, \cite{PVpoly}). In particular,
it ensures that in such domains a biharmonic function satisfies
the estimate
\begin{equation}\label{eq0.1}
\|\nabla u\|_{L^\infty(\Omega)}\leq C \|\nabla
u\|_{L^\infty(\po)}.
\end{equation}

\noindent Direct analogues of this principle for higher order
equations in general domains are unknown  (see Problem~4.3, p.275,
in J.\,Ne\v cas's book \cite{Necas}). Not only the increase of the
order leads to the failure of the methods which work for the
second order equations, but the properties of the solutions
themselves become more involved.

To be more specific,  let $\Omega\subset\RR^n$ be a bounded domain
and consider the boundary value problem
\begin{equation}\label{eq1.1}
\Delta^2 u= f \,\,{\mbox{in}}\,\,\Omega, \quad u\in \ring
W_2^2(\Omega),
\end{equation}

\noindent where the Sobolev space $\ring W_2^2(\Omega)$ is a
completion of $C_0^\infty(\Omega)$ in the norm $\|u\|_{\ring
W_2^2(\Omega)}=\|\Delta u\|_{L^2(\Omega)}$ and $f$ is a reasonably
nice function. Motivated by (\ref{eq0.1}), we ask if the gradient
of a solution to problem (\ref{eq1.1}) is bounded in an arbitrary
domain $\Omega\subset \RR^n$. It turns out that this property may
fail when $n\geq 4$ (see the counterexamples built in \cite{MR}
and \cite{PVLp}). In dimension three the boundedness of the
gradient of a solution has been an open problem.

The absence of any information about the geometry of the domain
puts this question beyond the scope of applicability of the
previously devised methods -- the aforementioned results regarding
the maximum principle heavily relied on specific assumptions on
$\Omega$. In the present paper we develop a new set of techniques
which allows to establish {\it the boundedness of the gradient of
the solution to (\ref{eq1.1}) under no restrictions on the
underlying domain}. Moreover, we prove the following:

\begin{theorem}\label{t1.1}
Let $\Omega$ be an arbitrary bounded domain in $\RR^3$ and let $G$
be Green's function for the biharmonic equation. Then
\begin{eqnarray}\label{eq1.3}
&&\qquad\qquad |\nabla_x\nabla_y G(x,y)|\leq C|x-y|^{-1},\qquad
x,y\in\Omega,\\[4pt]\label{eq1.3.1}
&&|\nabla_xG(x,y)|\leq C\quad\mbox{and}\quad|\nabla_yG(x,y)|\leq
C,\qquad x,y\in\Omega,
\end{eqnarray}

\noindent where $C$ is an absolute constant.
\end{theorem}

The boundedness of the gradient of a solution to the biharmonic
equation is a sharp property in the sense that the function $u$
satisfying (\ref{eq1.1}) generally does not exhibit more regularity.
Indeed, let $\Omega$ be the three-dimensional punctured unit ball
$B_1\setminus\{O\}$, where $B_r=\{x\in\RR^3:\,|x|<r\}$, and consider
a function $\eta\in C_0^\infty(B_{1/2})$ such that $\eta=1$ on
$B_{1/4}$. Let
\begin{equation}\label{eq1.4}
u(x):=\eta(x)|x|,\qquad x\in B_1\setminus\{O\}.
\end{equation}

\noindent Obviously, $u\in \ring W_2^2(\Omega)$ and $\Delta^2 u\in
C_0^\infty (\Omega)$. While $\nabla u$ is bounded, it is not
continuous at the origin. Therefore, the {\it continuity} of the
gradient {\it does not hold} in general and must depend on some
delicate properties of the domain.

Even in the case of the Laplacian the issue of continuity is
subtle. It has been resolved in 1924, when Wiener gave his famous
criterion for the regularity of a boundary point \cite{Wiener}.
Needless to say, Wiener's result strongly influenced the
development of partial differential equations, the theory of
function spaces and probability. Over the years it has been
extended to a variety of second order elliptic and parabolic
equations (\cite{LSW}, \cite{FJK}, \cite{FGL}, \cite{DMM},
\cite{MZ}, \cite{AHe}, \cite{TW}, \cite{La}, \cite{EG}; see also
the review papers \cite{MWiener}, \cite{AWiener}). However, the
case of higher-order operators is far from being well-understood.

Let us recall the original Wiener's criterion. Roughly speaking, it
states that a point $O\in\po$ is regular (i.e. every solution to the
Dirichlet problem for the Laplacian, with continuous data, is
continuous at $O$) if and only if the complement of the domain near
the point $O$, measured in terms of the Wiener (harmonic) capacity,
is sufficiently massive. More specifically, the harmonic capacity of
a compactum $K\subset\RR^n$ can be defined as
\begin{equation}\label{eq1.5}
{\rm cap}\,(K):=\inf\Bigl\{\|\nabla u\|_{L^2(\RR^n)}^2:\,\,u\in
C_0^\infty(\RR^n),\,\,u=1\mbox{ in a neighborhood of }K\Bigr\},
\end{equation}

\noindent where $n\geq 3$, and the regularity of the point $O$ is
equivalent to the condition
\begin{equation}\label{eq1.6}
\sum_{j=0}^\infty2^{j(n-2)}\, {\rm
cap}\,(\overline{B_{2^{-j}}}\setminus\Omega)=+\infty,
\end{equation}
where $B_{2^{-j}}$ is the ball of radius $2^{-j}$ centered at the origin. An appropriate version of this condition is also available in dimension $n=2$.

Recently, some progress has been made in the study of the
continuity of solutions for a certain family of higher order
elliptic equations in \cite{M2} (see also \cite{M0}, \cite{M1}).
In particular, these developments extend (\ref{eq1.6}) to the
context of the biharmonic equation in dimensions 4, 5, 6 and 7,
with the potential-theoretic capacity of order four in place of
(\ref{eq1.5}). In the present paper we pursue a different goal --
to obtain an analogue of the Wiener's test governing the {\it
gradient} of the solution.

Turning to this issue, we start with a suitable notion of
capacity. Let $\Pi$ denote the space of functions
\begin{equation}\label{eq1.7}
{\textstyle{P(x)=b_0+b_1\,\frac{x_1}{|x|}+b_2\,\frac{x_2}{|x|}+b_3\,
\frac{x_3}{|x|}}}, \qquad x\in\RR^3\setminus \{O\},\qquad
b_i\in\RR,\qquad i=0,1,2,3,
\end{equation}

\noindent  and $\Pi_1:=\{P\in\Pi:\,\|P\|_{\Pi}=1\}$. Then, given a
compactum $K\subset \RR^3\setminus\{0\}$ and $P\in\Pi_1$, let
\begin{equation}\label{eq1.8}
{\rm Cap}_P\,(K):=\inf\Bigl\{\|\Delta u\|_{L^2(\RR^3)}^2:\,\,u\in
\ring W_2^2(\RR^3\setminus\{0\}),\,\,u=P\mbox{ in a neighborhood
of }K\Bigr\}.
\end{equation}

\noindent This capacity first appeared in \cite{MT}, in the upper
estimates on $\sup_r ( \frac{1}{r^3}\int_{B_r}|\nabla
u(x)|^6\,dx)^{1/6}$ for a solution of (\ref{eq1.1}).

We say that a point $O\in\po$ is 1-regular if for every $f\in
C_0^\infty(\Omega)$ the solution $u$ to (\ref{eq1.1}) is
continuously differentiable at $O$, i.e. $\nabla u(x)\to 0$ as
$x\to O$; and $O$ is 1-irregular otherwise. Our main result
concerning 1-regularity is the following.

\begin{theorem}\label{t1.2} Let $\Omega$ be an open set in $\RR^3$.
If  for some $a\geq 2$
\begin{equation}\label{eq1.9}
\sum_{j=0}^\infty a^{-j} \,\inf_{P\in\Pi_1}{{\rm
Cap}_P\,(\overline{C_{a^{-j},a^{-j+1}}}\setminus\Omega)} =+\infty
\end{equation}

\noindent then the point $O$ is 1-regular.

Conversely, if the point $O\in\po$ is 1-regular then for every
$a\geq 2$
\begin{equation}\label{eq1.10}
\inf_{P\in\Pi_1}\sum_{j=0}^\infty a^{-j} \,{{\rm
Cap}_P\,(\overline{C_{a^{-j},a^{-j+2}}}\setminus\Omega)} =+\infty.
\end{equation}

\noindent Here and throughout the paper, $C_{s,bs}$ is the annulus $\{x\in\RR^3:\,s<|x|<bs\}$, $s>0$, $b>1$.
\end{theorem}

In \S{9} we further discuss the discrepancy between conditions
(\ref{eq1.9}) and (\ref{eq1.10}) and show by counterexample that
(\ref{eq1.9}) is not always necessary for 1-regularity.

To the best of our knowledge, Theorem~\ref{t1.2} is the first
Wiener-type result addressing the continuity of the derivatives of
a solution. It is accompanied by corresponding estimates, in
particular, we prove the following refinement of (\ref{eq1.3}). Let $\Omega$ be a bounded domain in $\RR^3$,
$O\in\po$. Fix some $a\geq 4$ and let $c_a:=1/(32a^4)$. Then for
$x,y\in\Omega$
\begin{eqnarray}
\nonumber &&|\nabla_x\nabla_y G(x,y)|\\[4pt]\nonumber
&& \quad \leq  \left\{\begin{array}{l} \frac{C}{|x-y|}
\times\exp \left( -c\sum_{j=2}^{l_{yx}}(|y|a^{2j}) {\rm
Cap}\,(\overline{C_{32|y|a^{2(j-1)},32|y|a^{2j}}}\setminus\Omega)\right),\\[8pt]\nonumber
 \quad\quad\mbox{if  $|y|\leq c_a |x|$ and \,$l_{yx}\geq 2$, $l_{yx}\in\NN$, is such that $|x|\geq 32 a^{2l_{yx}}|y|$},\\[8pt]
\frac{C}{|x-y|} \times\exp \left( -c\sum_{j=2}^{l_{xy}}(|x|a^{2j}) {\rm
Cap}\,(\overline{C_{32|x|a^{2(j-1)},32|x|a^{2j}}}\setminus\Omega)\right),\\[8pt]\nonumber
 \quad\quad\mbox{if  $|x|\leq c_a |y|$ and \,$l_{xy}\geq 2$, $l\in\NN$, is such that $|y|\geq 32 a^{2l_{xy}}|x|$},\\[8pt]
\frac{C}{|x-y|}, \qquad \mbox{if}\quad c_a |y|\leq |x|\leq
c_a^{-1}|y|.
\end{array}
\right.
\end{eqnarray}

It has to be noted that Theorem~\ref{t1.2} brings up a peculiar
role of circular cones and planes for 1-regularity of a boundary
point. For example, if the complement of $\Omega$ is a compactum
located on the circular cone (or plane) given by
$\{x\in\RR^3\setminus\{0\}:\,b_0 |x|+b_1x_1+b_2x_2+b_3x_3=0\}$
such that the harmonic capacity ${\rm
cap}\,(\RR^3\setminus\Omega)=0$, then ${\rm
Cap}_P(\RR^3\setminus\Omega)=0$ for $P$ associated to the same
$b_i$'s. Hence, by Theorem~\ref{t1.2}, the point $O$ is not
1-regular.

Another surprising effect, strikingly different from the classical
theory, is that for some domains 1-irregularity turns out to be
unstable under affine transformations of coordinates.

In conclusion, we provide some examples further illustrating the
geometric nature of conditions (\ref{eq1.9})--(\ref{eq1.10}). Among
them is the model case when $\Omega$ has an inner cusp, i.e. in a
neighborhood of the origin $\Omega=\{(r,\theta,\phi): \,0<r<c,
\,h(r)<\theta\leq \pi,\, 0\leq\phi<2\pi\}$, where $h$ is a
non-decreasing function such that $h(br)\leq h(r)$ for some $b>1$.
For such a domain Theorem~\ref{t1.2} yields the following criterion:
\begin{equation}\label{eq1.11}
\mbox{the point $O$ is 1-regular \quad if and only if }\quad
\int_{0}^{1} s^{-1}h(s)^2\,ds=+\infty.
\end{equation}

\noindent Some other geometrical examples can be found in the body
of the paper.

\section{Integral identity and global estimate}
\setcounter{equation}{0}

Let us start with a few remarks about the notation.

Let $(r,\omega)$ be spherical coordinates in $\RR^3$, i.e. $r=|x|\in
(0,\infty)$ and $\omega=x/|x|$ is a point of the unit sphere $S^2$.
Occasionally we will write the spherical coordinates as
$(r,\theta,\phi)$, where $\theta\in[0,\pi]$ stands for the
colatitude and $\phi\in [0,2\pi)$ is the longitudinal coordinate,
i.e.
\begin{equation}\label{eq2.1}
\omega=x/|x|=(\sin
\theta\cos\phi,\,\sin\theta\sin\phi,\,\cos\theta).
\end{equation}

\noindent Now let $t=\log r^{-1}$. Then by  $\kappa$ and
$\varkappa$ we denote the mappings
\begin{equation}\label{eq2.2}
\RR^3 \ni x\,\stackrel{\kappa}{\longrightarrow}
\,(r,\phi,\theta)\in [0,\infty)\times [0,2\pi) \times
[0,\pi];\qquad\quad\RR^3 \ni
x\,\stackrel{\varkappa}{\longrightarrow} \,(t,\omega)\in\RR\times
S^{2}.
\end{equation}

\noindent The symbols $\delta_\omega$ and $\nabla_\omega$ refer,
respectively, to the Laplace-Beltrami operator and the gradient on
$S^2$.

For any domain $\Omega\subset \RR^3$ a function  $u\in
C_0^\infty(\Omega)$ can be extended by zero to $\RR^3$ and we will
write $u\in C_0^\infty(\RR^3)$ whenever convenient. Similarly, the
functions in $\ring W_2^2(\Omega)$ will be extended by zero and
treated as functions on $\RR^3$ without further comments.

By $C$, $c$, $C_i$ and $c_i$, $i\in\NN$, we generally denote some
constants whose exact values are of no importance. Also, we write
$A\approx B$, if $C^{-1}\,A\leq B\leq C\,A$ for some $C>0$.

The first result is

\begin{lemma}\label{l2.1}
Let $\Omega$ be an open set in $\RR^3$, $u\in C_0^\infty(\Omega)$
and $v=e^t(u\circ \varkappa^{-1})$. Then
\begin{eqnarray}\label{eq2.3}
&&\hskip -1cm \int_{\RR^3}\Delta
u(x)\,\Delta\bigg(u(x)|x|^{-1}\,{\cal G}(\log |x|^{-1})\bigg)\,dx\nonumber\\[4pt]
&&\hskip -1cm\quad\quad  =\int_{\RR}\int_{S^{2}}\Bigl[
(\delta_\omega v)^2{\cal G}+2(\partial_t\nabla_\omega v)^2{\cal
G}+ (\partial_t^2v)^2{\cal G}   -(\nabla_\omega v)^2
\Bigl(\partial_t^2
{\cal G}+\partial_t{\cal G}+2{\cal G}\Bigr) \nonumber\\[4pt]
&&\hskip -1cm\quad\quad -(\partial_t v)^2\Bigl(2\partial_t^2{\cal
G}+3\partial_t {\cal G}-{\cal G}\Bigr) +\frac
12\,v^2\Bigl(\partial_t^4{\cal G}+2\partial_t^3{\cal G}
-\partial_t^2{\cal G}-2\partial_t{\cal G}\Bigr)\Bigr]\,d\omega dt,
\end{eqnarray}

\noindent for every  function ${\cal G}$ on $\RR$ such that both
sides of {\rm (\ref{eq2.3})} are well-defined.
\end{lemma}

\vskip 0.08in

\bp In the system of coordinates $(t,\omega)$ the $3$-dimensional
Laplacian can be written as
\begin{equation}\label{eq2.4}
\Delta=e^{2t}\Lambda(\partial_t,\delta_\omega),
\quad\mbox{where}\quad
\Lambda(\partial_t,\delta_\omega)=\partial_t^2-\partial_t+\delta_\omega.
\end{equation}

\noindent Then passing to the coordinates $(t,\omega)$, we have
\begin{eqnarray}\label{eq2.5}
&&\hskip -1.3cm \int_{\RR^3}\Delta
u(x)\,\Delta\bigg(u(x)|x|^{-1}{\cal G}(\log |x|^{-1})\bigg)\,dx =
\int_{\RR}\int_{S^2} \Lambda(\partial_t-1,\delta_\omega)
v\,\,\Lambda(\partial_t,\delta_\omega)(v{\cal G})\,d\omega dt\nonumber\\[4pt]
&&\hskip -1.3cm\quad = \int_{\RR}\int_{S^2}
\left(\partial_t^2v-3\partial_tv+2v+\delta_\omega v\right)
\,\left(\partial_t^2(v{\cal G})-\partial_t(v{\cal G})+{\cal
G}\,\delta_\omega
v\right)\,d\omega dt\nonumber\\[4pt]
&&\hskip -1.3cm\quad = \int_{\RR}\int_{S^{2}}
\left(\partial_t^2v-3\partial_tv+2v+\delta_\omega v\right)\nonumber\\[4pt]
&&\quad\quad \times\left({\cal G}\,\delta_\omega v+{\cal
G}\,\partial_t^2v+(2\partial_t {\cal G}-{\cal
G})\,\partial_tv+(\partial_t^2{\cal G}-\partial_t {\cal
G})\,v\right)\,d\omega
dt\nonumber\\[4pt]
&&\hskip -1.3cm\quad= \int_{\RR}\int_{S^{2}} \Bigl(
\left((\delta_\omega v)^2+2\,\delta_\omega
v\partial_t^2 v+\,(\partial_t^2v)^2\right){\cal G} \nonumber\\[4pt]
&&\hskip -1.3cm\quad\quad+\left(v\delta_\omega v+v\partial_t^2
v\right)\,\left(\partial_t^2 {\cal G}-\partial_t{\cal G}+2{\cal
G}\right)+\left(\delta_\omega v\partial_t v+\partial_t^2
v\partial_t v\right)\left(2\partial_t {\cal G}-4{\cal G}\right)
\nonumber\\[4pt]
&&\hskip -1.3cm\quad\quad+(\partial_tv)^2\,\left(-6\partial_t
{\cal G}+3{\cal G}\right) +v\partial_t v\left(-3\partial_t^2{\cal
G}+7\partial_t{\cal G}-2{\cal
G}\right)+v^2\left(2\partial_t^2{\cal G}-2\partial_t{\cal
G}\right) \Bigr)\,d\omega dt.
\end{eqnarray}

\noindent This, in turn, is equal to
\begin{eqnarray}\label{eq2.6}
&&\hskip -0.7cm  \int_{\RR}\int_{S^{2}} \Bigl( {\cal
G}\,(\delta_\omega v)^2-2{\cal G}\,\delta_\omega\partial_t
v\,\partial_t v+{\cal G}\,(\partial_t^2v)^2\nonumber\\[4pt]
&&\hskip -0.7cm\quad+(\nabla_\omega v)^2\,\left(-\partial_t^2{\cal
G}-(\partial_t^2 {\cal G}-\partial_t{\cal G}+2{\cal
G})+(\partial_t^2 {\cal G}-2\partial_t{\cal G})\right)
\nonumber\\[4pt]
&&\hskip -0.7cm\quad +(\partial_t v)^2\left(-(\partial_t^2 {\cal
G}-\partial_t{\cal G}+2{\cal G})+(-\partial_t^2{\cal
G}+2\partial_t{\cal G})+(-6\partial_t {\cal G}+3{\cal G})\right)
\nonumber\\[4pt]
&&\hskip -0.7cm\quad+v\partial_t v\left(-(\partial_t^3 {\cal
G}-\partial_t^2{\cal G}+2\partial_t{\cal G})+(-3\partial_t^2{\cal
G}+7\partial_t{\cal G}-2{\cal G})
\right)\nonumber\\[4pt]
&&\hskip -0.7cm\quad+v^2\left(2\partial_t^2{\cal
G}-2\partial_t{\cal G}\right) \Bigr)\,d\omega dt,
\end{eqnarray}

\noindent and integrating by parts once again we obtain
(\ref{eq2.3}). \ep

In order to single out the term with $v^2$ in (\ref{eq2.3}) we
shall need the following auxiliary result.

\begin{lemma}\label{l2.2} Consider the equation
\begin{equation}\nonumber
\frac{d^4g}{dt^4}+2\frac{d^3g}{dt^3}
-\frac{d^2g}{dt^2}-2\frac{dg}{dt} =\delta,  \label{eq2.7}
\end{equation}

\noindent where $\delta$ stands for the Dirac delta function. A
unique solution to {\rm (\ref{eq2.7})} which is bounded and vanishes
at $+\infty$ is given by
\begin{equation}\label{eq2.8}
g(t)=-\frac 16\left\{\begin{array}{l}
e^t-3,\qquad \qquad\qquad t<0,\\[4pt]
e^{ -2t} -3\,e^{-t},\qquad\,\,\,\quad t>0.\\[4pt]
\end{array}
\right.
\end{equation}

\end{lemma}

\bp Since the equation (\ref{eq2.7}) is equivalent to
\begin{equation}\label{eq2.9}
\frac{d}{dt}\left(\frac{d}{dt}+2\right)\left(\frac{d}{dt}+1\right)
\left(\frac{d}{dt}-1\right)g=\delta,
\end{equation}

\noindent a bounded solution of (\ref{eq2.7}) vanishing at
$+\infty$ must have the form
\begin{equation}\label{eq2.10}
g(t)=\left\{\begin{array}{l}
a\, e^t+b,\qquad \qquad\qquad t<0,\\[4pt]
c\, e^{ -2t} +d\,e^{-t},\qquad\quad\,\,\, t>0,\\[4pt]
\end{array}
\right.
\end{equation}

\noindent for some constants $a,b,c,d$. Once this is established,
we find the system of coefficients so that $\partial_t^kg$ is
continuous for $k=0,1,2$ and $\lim_{t\to 0^+}\partial_t^3
g(t)-\lim_{t\to 0^-}\partial_t^3 g(t)=1$.
 \ep

With Lemma~\ref{l2.2} at hand, a suitable choice of the function
${\cal G}$ yields the positivity of the left-hand side of
(\ref{eq2.3}), one of the cornerstones of this paper. The details
are as follows.

\begin{lemma}\label{l2.3}
Let $\Omega$ be a bounded domain in $\RR^3$, $O\in
\RR^3\setminus\Omega$, $u\in C_0^\infty(\Omega)$ and $v=e^t(u\circ
\varkappa^{-1})$. Then for every $\xi\in\Omega$ and $\tau=\log
|\xi|^{-1}$ we have

\begin{equation}\label{eq2.11}
\frac 12\int_{S^{n-1}}v^2(\tau,\omega)\,d\omega \leq
\int_{\RR^n}\Delta u(x)\,\Delta\bigg(u(x)|x|^{-1}g(\log
(|\xi|/|x|))\bigg)\,dx,
\end{equation}

\noindent where $g$ is given by {\rm (\ref{eq2.8})}.

\end{lemma}

\bp Representing $v$ as a series of spherical harmonics and noting
that the eigenvalues
 of the Laplace-Beltrami operator  on the unit sphere are
$k(k+1)$, $k=0,1,...$, we arrive at the inequality
\begin{equation}\label{eq2.12}
\int_{S^2}|\delta_\omega v|^2\,d\omega\geq 2
\int_{S^2}|\nabla_\omega v|^2\,d\omega.
\end{equation}

Now, let us take ${\cal G}(t)=g(t-\tau)$, $t\in\RR$. Since $g\geq
0$, the combination of Lemma \ref{l2.2}, (\ref{eq2.3}) and
(\ref{eq2.12}) allows one to obtain the estimate
\begin{eqnarray}\label{eq2.13}
&&\hskip -1cm \int_{\RR^n}\Delta
u(x)\,\Delta\bigg(u(x)|x|^{-1}g(\log
(|\xi|/|x|))\bigg)\,dx\nonumber\\[4pt]
&&\hskip -1cm\quad \geq \int_{\RR}\int_{S^{n-1}}\Bigl[
 -(\nabla_\omega v(t,\omega))^2 \Bigl(\partial_t^2
g(t-\tau)+\partial_tg(t-\tau)\Bigr)\nonumber\\[4pt]
&&\hskip -1cm\quad  -(\partial_t
v(t,\omega))^2\Bigl(2\partial_t^2g(t-\tau)+3\partial_t
g(t-\tau)-g(t-\tau)\Bigr)\Bigr]\,d\omega dt +\frac
12\int_{S^{n-1}}v^2(\tau,\omega)\,d\omega.
\end{eqnarray}

\noindent Thus, the matters are reduced to showing that
\begin{equation}\label{eq2.14}
\partial_t^2
g+\partial_tg\leq 0\qquad\mbox{and}\qquad
2\partial_t^2g+3\partial_t g-g\leq 0.
\end{equation}

\noindent Indeed, we compute
\begin{equation}\label{eq2.15}
\partial_tg(t)=-\frac 16\left\{\begin{array}{l}
e^t,\qquad \qquad\qquad \quad \,t<0,\\[4pt]
-2e^{ -2t} +3\,e^{-t},\qquad\, t>0,
\end{array}
\right.
\end{equation}

\noindent and
\begin{equation}\label{eq2.16}
\partial_t^2g(t)=-\frac 16\left\{\begin{array}{l}
e^t,\qquad \qquad\qquad \quad t<0,\\[4pt]
4e^{ -2t} -3\,e^{-t},\,\,\qquad t>0,
\end{array}
\right.
\end{equation}

\noindent which gives
\begin{equation}\label{eq2.17}
\partial_t^2
g(t)+\partial_tg(t)=-\frac 13\left\{\begin{array}{l}
e^t,\qquad  \quad t<0,\\[4pt]
e^{ -2t},\qquad t>0,
\end{array}
\right.
\end{equation}

\noindent and
\begin{equation}\label{eq2.18}
2\partial_t^2g(t)+3\partial_t g(t)-g(t)=-\frac
16\left\{\begin{array}{l}
4e^t+3,\qquad  \quad t<0,\\[4pt]
e^{ -2t}+6e^{-t},\quad\,\, t>0.
\end{array}
\right.
\end{equation}

\noindent Clearly, both functions (\ref{eq2.17}), (\ref{eq2.18}) are
non-positive. The result follows from (\ref{eq2.13}).
  \ep

\section{Local energy and $L^2$ estimates}
\setcounter{equation}{0}

This section is devoted to estimates for a solution of the Dirichlet
problem near a boundary point, in particular, the proof of
Theorem~\ref{t1.1}. To set the stage, let us first record the
well-known result following from the energy estimate for solutions
of elliptic equations.

\begin{lemma}\label{l3.1} Let $\Omega$ be an arbitrary
domain in $\RR^3$, $Q\in\RR^3\setminus\Omega$ and $R>0$. Suppose
\begin{equation}\label{eq3.1}
\Delta^2 u=f \,\,{\mbox{in}}\,\,\Omega, \quad f\in
C_0^{\infty}(\Omega\setminus B_{4R}(Q)),\quad u\in \ring
W_2^2(\Omega).
\end{equation}

\noindent Then
\begin{equation}\label{eq3.2}
\int_{B_{\rho}(Q)\cap\Omega}|\nabla^2 u|^2\,dx
+\frac{1}{\rho^2}\int_{B_{\rho}(Q)\cap\Omega}|\nabla u|^2\,dx\leq
\frac{C}{\rho^4}\int_{C_{\rho,2\rho}(Q)\cap\Omega}|u|^2\,dx
\end{equation}

\noindent for every $\rho<2R$.
\end{lemma}

Here and throughout the paper $B_r(Q)$ and $S_r(Q)$ denote,
respectively, the ball and the sphere with radius $r$ centered at
$Q$ and $C_{r,R}(Q)=B_R(Q)\setminus\overline{B_r(Q)}$. When the
center is at the origin, we write $B_r$ in place of $B_r(O)$, and
similarly $S_r:=S_r(O)$ and $C_{r,R}:=C_{r,R}(O)$. Also, $\nabla^2u$
stands for a vector of all second derivatives of $u$.

We omit a standard proof of Lemma~\ref{l3.1} (see, e.g., \cite{ADN},
\cite{Shen}) and proceed to estimates for a biharmonic function
based upon the results in \S{2}.

\begin{proposition}\label{p3.2}
Let $\Omega$ be a bounded domain in $\RR^3$,
$Q\in\RR^3\setminus\Omega$, and $R>0$. Suppose
\begin{equation}\label{eq3.3}
\Delta^2 u=f \,\,{\mbox{in}}\,\,\Omega, \quad f\in
C_0^{\infty}(\Omega\setminus B_{4R}(Q)),\quad u\in \ring
W_2^2(\Omega).
\end{equation}

\noindent Then
\begin{equation}\label{eq3.4}
\frac{1}{\rho^{4}}\int_{S_{\rho}(Q)\cap\Omega}|u(x)|^2\,d\sigma_x
\leq \frac{C}{R^5} \int_{C_{R,4R}(Q)\cap\Omega} |u(x)|^2\,dx\quad
{\mbox{ for every}}\quad \rho<R,
\end{equation}

\noindent where $C$ is an absolute constant.
\end{proposition}

\bp For notational convenience we assume that $Q=O$. Let us
approximate $\Omega$ by a sequence of domains with smooth boundaries
$\left\{\Omega_n\right\}_{n=1}^\infty$ satisfying
\begin{equation}\label{eq3.5}
\bigcup_{n=1}^\infty \Omega_n=\Omega\quad\mbox{and}\quad
{\overline{\Omega}}_n\subset\Omega_{n+1} \quad \mbox{for every}
\quad n\in\NN.
\end{equation}

\noindent Choose $n_0\in\NN$ such that ${\rm supp}\,f\subset
\Omega_n$ for every $n\geq n_0$ and denote by $u_n$ a unique
solution of the Dirichlet problem
\begin{equation}\label{eq3.6}
\Delta^2 u_n = f\quad {\rm in} \quad \Omega_n, \quad  u_n\in \ring
W_2^2(\Omega_n), \quad n\geq n_0.
\end{equation}

\noindent The sequence  $\{u_n\}_{n=n_0}^\infty$ converges to $u$ in
$\ring W_2^2(\Omega)$ (see, e.g., \cite{Necas}, \S{6.6}).

Next,  take some $\eta\in C_0^\infty(B_{2R})$ such that
\begin{equation}\label{eq3.7}
0\leq\eta\leq 1\,\, {\rm in}\,\,B_{2R}, \quad \eta=1\,\, {\rm
in}\,\,B_R\quad {\rm and} \quad |\nabla^k \eta|\leq CR^{-k}, \quad
k\leq 4.
\end{equation}

\noindent  Also, fix $\tau=\log\rho^{-1}$ and let $g$ be the
function defined in (\ref{eq2.8}).

Consider the difference
\begin{eqnarray}\label{eq3.8}
&&\int_{\RR^3}\Delta \bigg(\eta(x)u_n(x)\bigg)\,
\Delta\bigg(\eta(x)u_n(x)|x|^{-1}g(\log (\rho/|x|))\bigg)
\,dx\nonumber\\[4pt]
&&\quad  -\int_{\RR^3}\Delta
u_n(x)\,\Delta\bigg(u_n(x)|x|^{-1}g(\log
(\rho/|x|))\eta^2(x)\bigg)\,dx.
\end{eqnarray}

\noindent One can view this expression as
\begin{equation}\label{eq3.9}
\int_{\RR^3}\bigg([\Delta^2,\eta] u_n(x)\bigg) \bigg(\eta(x)u_n(x)|x|^{-1}g(\log (\rho/|x|))\bigg)\,dx,
\end{equation}

\noindent where the integral is understood in the sense of pairing
between $\ring W_2^2(\Omega_n)$ and its dual. Evidently, the
support of the integrand is a subset of ${\rm
supp}\,\nabla\eta\subset C_{R,2R}$, and therefore, the difference
in (\ref{eq3.8}) is bounded by
\begin{equation}\label{eq3.10}
C\sum_{k=0}^2 \frac{1}{R^{5-2k}}\int_{C_{R,2R}} |\nabla^k
u_n(x)|^2\,dx.
\end{equation}

Since $u_n$ is biharmonic in $\Omega_n\cap B_{4R}$ and $\eta$ is
supported in $B_{2R}$, the second term in (\ref{eq3.8}) is equal to
zero. Turning to the first term, we shall employ Lemma~\ref{l2.3}
with $u=\eta\,u_n$. The result of the Lemma holds for such a choice
of $u$. This can be seen directly by inspection of the argument or
one can approximate each $u_n$ by a sequence of
$C_0^\infty(\Omega_n)$ functions in $\ring W_2^2(\Omega_n)$ and then
take a limit using that $O\not\in\overline\Omega_n$. Then
(\ref{eq3.8}) is bounded from below by
\begin{equation}\label{eq3.11}
\frac{C}{\rho^{4}}\int_{S_{\rho}}|\eta(x)u_n(x)|^2\,d\sigma_x.
\end{equation}

\noindent Hence, for every $\rho<R$
\begin{equation}\label{eq3.12}
\frac{1}{\rho^{4}}\int_{S_{\rho}}|u_n(x)|^2\,d\sigma_x \leq C
\sum_{k=0}^2 \frac{1}{R^{5-2k}}\int_{C_{R,2R}} |\nabla^k
u_n(x)|^2\,dx.
\end{equation}

\noindent Now the proof can be finished applying Lemma~\ref{l3.1}
and taking the limit as $n\to \infty$.  \ep

Now we show that (\ref{eq3.4}) yields a uniform  pointwise estimate
for $\nabla u$.

\begin{corollary}\label{c3.3}
Let $\Omega$ be a bounded domain in $\RR^3$,
$Q\in\RR^3\setminus\Omega$, $R>0$ and
\begin{equation}\label{eq3.13}
\Delta^2 u=f \,\,{\mbox{in}}\,\,\Omega, \quad f\in
C_0^{\infty}(\Omega\setminus B_{4R}(Q)),\quad u\in \ring
W_2^2(\Omega).
\end{equation}

\noindent Then for every $x\in B_{R/4}(Q)\cap\Omega$
\begin{equation}\label{eq3.14}
|\nabla u(x)|^2\leq\frac{C}{R^{5}}\int_{C_{R/4,4R}(Q)\cap\Omega}
|u(y)|^2\,dy,
\end{equation}

\noindent and
\begin{equation}\label{eq3.15}
|u(x)|^2\leq C \frac{|x-Q|^2}{R^{5}}\int_{C_{R/4,4R}(Q)\cap\Omega}
|u(y)|^2\,dy.
\end{equation}

In particular, for every bounded domain $\Omega\subset\RR^3$ the
solution to the boundary value problem {\rm (\ref{eq3.13})}
satisfies
\begin{equation}\label{eq3.16}
|\nabla u|\in L^\infty(\Omega).
\end{equation}
\end{corollary}

\bp By an interior estimate for solutions of the elliptic equations
(see \cite{FJohn}, pp. 153-155)
\begin{equation}\label{eq3.17}
|\nabla u(x)|^2\leq \frac{C}{d(x)^3} \int_{B_{d(x)/2}(x)} |\nabla
u(y)|^2\,dy,
\end{equation}

\noindent where $d(x)$ denotes the distance from $x$ to $\po$. Let
$x_0$ be a point on the boundary of $\Omega$ such that
$d(x)=|x-x_0|$. Since $x\in B_{R/4}(Q)\cap\Omega$ and
$Q\in\RR^3\setminus\Omega$, we have $x\in B_{R/4}(x_0)$, and
therefore
\begin{equation}\label{eq3.18}
\frac{1}{d(x)^3}\int_{B_{d(x)/2}(x)} |\nabla u(y)|^2\,dy\leq
\frac{C}{d(x)^5} \int_{B_{ 2d(x)}(x_0)} |u(y)|^2\,dy \leq
\frac{C}{R^{5}}\int_{C_{3R/4,3R}(x_0)} |u(y)|^2\,dy,
\end{equation}

\noindent using Lemma~\ref{l3.1} for the first estimate and
(\ref{eq3.4}) for the second one. Indeed, $d(x)\leq R/4$ and
therefore, $2d(x)< 3R/4$. On the other hand, $u$ is biharmonic in
$B_{4R}(Q)\cap\Omega$ and
\begin{equation}\label{eq3.19}
|Q-x_0|\leq |Q-x|+|x-x_0|\leq R/2.
\end{equation}

\noindent Hence, $u$ is biharmonic in $B_{3R}(x_0)\cap\Omega$ and
Proposition~\ref{p3.2} holds with $x_0$ in place of $Q$, $3R/4$ in
place of $R$ and $\rho=2d(x)$. Furthermore, (\ref{eq3.19}) yields
\begin{equation}\label{eq3.20}
C_{3R/4,3R}(x_0)\subset C_{R/4,4R}(Q),
\end{equation}

\noindent and that finishes the argument for (\ref{eq3.14}).

To prove (\ref{eq3.15}), we start with the estimate
\begin{equation}\label{eq3.21}
|u(x)|^2\leq \frac{C}{d(x)^3} \int_{B_{d(x)/2}(x)} |u(y)|^2\,dy,
\end{equation}

\noindent and then proceed using (\ref{eq3.4}), much as in
(\ref{eq3.18})--(\ref{eq3.20}).  \ep

Using the Kelvin transform for biharmonic functions, an estimate on
a biharmonic function near the origin can be translated into an
estimate at infinity. In particular, Proposition~\ref{p3.2} and
Corollary~\ref{c3.3} lead to the following result.

\begin{proposition}\label{p3.4} Let $\Omega$ be a bounded domain in $\RR^3$,
$Q\in\RR^3\setminus\Omega$, $r>0$ and assume that
\begin{equation}\label{eq3.22}
\Delta^2 u=f \,\,{\mbox{in}}\,\,\Omega, \quad f\in
C_0^{\infty}(B_{r/4}(Q)\cap\Omega),\quad u\in \ring W_2^2(\Omega).
\end{equation}

\noindent Then
\begin{equation}\label{eq3.23}
\frac{1}{\rho^{2}}\int_{S_{\rho}(Q)\cap\Omega}|u(x)|^2\,d\sigma_x
\leq \frac{C}{r^{3}} \int_{ C_{r/4,r}(Q)\cap\Omega} |u(x)|^2\,dx,
\end{equation}

\noindent for any $\rho>r$.

Furthermore, for any $x\in\Omega\setminus B_{4r}(Q)$
\begin{equation}\label{eq3.24}
|\nabla u(x)|^2\leq
\frac{C}{|x-Q|^2\,r^3}\int_{C_{r/4,4r}(Q)\cap\Omega} |u(y)|^2\,dy,
\end{equation}

\noindent and
\begin{equation}\label{eq3.25}
|u(x)|^2\leq \frac{C}{r^3}\int_{C_{r/4,4r}(Q)\cap\Omega}
|u(y)|^2\,dy.
\end{equation}
\end{proposition}

\bp As before, it is enough to consider the case $Q=O$. Retain the
approximation of $\Omega$ with the sequence of smooth domains
$\Omega_n$ satisfying (\ref{eq3.5}) and define $u_n$ according to
(\ref{eq3.6}). We denote by ${\cal I}$ the inversion $x\mapsto
y=x/|x|^2$ and by $U_n$ the Kelvin transform of $u_n$,
\begin{equation}\label{eq3.26}
U_n(y):=|y|\, u_n(y/|y|^2), \quad y\in {\cal I}(\Omega_n).
\end{equation}

\noindent Then
\begin{equation}\label{eq3.27}
\Delta^2 U_n(y)=|y|^{-7} (\Delta^2 u_n)(y/|y|^2),
\end{equation}

\noindent and therefore, $U_n$ is biharmonic in ${\cal
I}(\Omega_n)\cap B_{4/r}$. Moreover, (\ref{eq3.27}) implies that
\begin{equation}\label{eq3.28}
\int_{{\cal I}(\Omega_n)} |\Delta U_n(y)|^2\, dy= \int_{\Omega_n}
|\Delta u_n(x)|^2\, dx,
\end{equation}

\noindent so that
\begin{equation}\label{eq3.29}
U_n\in \ring W_2^2({\cal I}(\Omega_n))
\quad\Longleftrightarrow\quad u_n\in \ring W_2^2(\Omega_n).
\end{equation}

\noindent Observe also that $\Omega_n$ is a bounded domain with
$O\not\in \overline{\Omega}_n$, hence, so is ${\cal I}(\Omega_n)$
and $O\not\in \overline{{\cal I}(\Omega_n)}$.

Following Proposition~\ref{p3.2}, we show that
\begin{equation}\label{eq3.30}
\rho^4\int_{S_{1/\rho}}|U_n(y)|^2\,d\sigma_y \leq C\,
r^5\int_{C_{1/r,4/r}} |U_n(y)|^2\,dy,
\end{equation}

\noindent which after the substitution (\ref{eq3.26}) and the
change of coordinates yields
\begin{equation}\label{eq3.31}
\frac{1}{\rho^{2}}\int_{S_{\rho}}|u_n(x)|^2\,d\sigma_x \leq
\frac{C}{r^{3}} \int_{ C_{r/4,r}} |u_n(x)|^2\,dx.
\end{equation}

Turning to the pointwise estimates (\ref{eq3.24})--(\ref{eq3.25}),
let us fix some $x\in\Omega\setminus B_{4r}(Q)$. Observe that
\begin{equation}\label{eq3.32}
|\nabla u_n(x)|\leq C|x|^{-1} \left|(\nabla
U_n)(x/|x|^2)\right|+\left|U_n(x/|x|^2)\right|,
\end{equation}

\noindent since $u_n(x)=|x|\, U_n(x/|x|^2)$. Therefore, combining
(\ref{eq3.32}) and Corollary~\ref{c3.3} applied to the function
$U_n$, we deduce that
\begin{equation}\label{eq3.33}
|\nabla u_n(x)|^2\leq C
\,\frac{r^5}{|x|^2}\int_{C_{1/(4r),4/r}}|U_n(z)|^2\,dz=
 \,\frac{C}{|x|^2\,r^3}\int_{C_{r/4,4r}}|u_n(z)|^2\,dz,
\end{equation}

\noindent and
\begin{equation}\label{eq3.34}
|u_n(x)|^2\leq C r^5\int_{C_{1/(4r),4/r}}|U_n(z)|^2\,dz=
\frac{C}{r^3}\int_{C_{r/4,4r}}|u_n(z)|^2\,dz.
\end{equation}

At this point, we can use the limiting procedure to complete the
argument. Indeed, since $u_n$ converges to $u$ in $\ring
W_2^2(\Omega)$, the integrals in (\ref{eq3.31}), (\ref{eq3.33})
and (\ref{eq3.34}) converge to the corresponding integrals with
$u_n$ replaced by $u$. Turning to $|\nabla u_n(x)|$, we observe
that both $u_n$ and $u$ are biharmonic in a neighborhood of $x$,
in particular, for sufficiently small $d$
\begin{equation}\label{eq3.35}
|\nabla (u_n(x)-u(x))|^2\leq
\,\frac{C}{d^5}\int_{B_{d/2}(x)}|u_n(z)-u(z)|^2\,dz.
\end{equation}

\noindent As $n\to\infty$, the integral on the right-hand side of
(\ref{eq3.35}) vanishes and therefore, $|\nabla u_n(x)|\to |\nabla
u(x)|$. Similar considerations apply to $u_n(x)$. \ep


\section{Estimates for Green's function}
\setcounter{equation}{0}

Let $\Omega$ be a bounded three-dimensional domain. As in the
introduction, we denote by $G(x,y)$, $x,y\in\Omega$, Green's
function for the biharmonic equation. In other words, for every
fixed $y\in\Omega$ the function $G(x,y)$ satisfies
\begin{equation}\label{eq4.1}
\Delta_x^2 G(x,y)=\delta(x-y), \qquad x\in\Omega,
\end{equation}

\noindent in the space $\ring W_2^2(\Omega)$. Here and throughout
the section $\Delta_x$ stands for the Laplacian in $x$ variable, and
similarly we use the notation $\Delta_y$, $\nabla_y$, $\nabla_x$ for
the Laplacian and gradient in $y$, and gradient in $x$,
respectively. As before, $d(x)$ is the distance from $x\in\Omega$ to
$\po$.

\begin{proposition}\label{p4.1} Let $\Omega\subset \RR^3$ be a bounded
domain. Then there exists an absolute constant $C$ such that for
every $x,y\in\Omega$
\begin{equation}\label{eq4.2}
\bigg|\nabla_x\nabla_y (G(x,y)-\Gamma(x-y))\bigg|\leq
\frac{C}{\max\{|x-y|,\,d(x),\,d(y)\}},
\end{equation}

\noindent where $\Gamma(x-y)=\frac{|x-y|}{8\pi}$ is the fundamental
solution for the bilaplacian.
\end{proposition}

\bp Let us start with some auxiliary calculations. Consider a
function $\eta$ such that
\begin{equation}\label{eq4.4}
\eta \in C_0^\infty(B_{1/2})\quad\mbox{and}\quad \eta=1
\quad\mbox{in}\quad B_{1/4},
\end{equation}

\noindent and define a vector-valued function ${\cal R}=({\cal
R}_1,{\cal R}_2,{\cal R}_3)$ by
\begin{equation}\label{eq4.5}
{\cal R}_j (x,y):=\frac{\partial}{\partial y_j}\,
G(x,y)-\eta\left(\frac{x-y}{d(y)}\right)\frac{\partial}{\partial
y_j}\, \Gamma (x-y),\qquad x,y\in\Omega,
\end{equation}

\noindent where $j=1,2,3$. Also, let us denote
\begin{equation}\label{eq4.6}
f_j(x,y):= \Delta_x^2 {\cal R}_j
(x,y)=-\left[\Delta_x^2,\eta\left(\frac{x-y}{d(y)}\right)
\right]\frac{\partial}{\partial y_j}\,\Gamma (x-y),\quad j=1,2,3.
\end{equation}

\noindent It is not hard to see that for every $j$
\begin{equation}\label{eq4.7}
f_j(\cdot,y)\in C_0^\infty(C_{d(y)/4,
d(y)/2}(y))\qquad\mbox{and}\qquad |f_j(x,y)|\leq Cd(y)^{-4}, \quad
x,y\in\Omega.
\end{equation}

\noindent Then for every fixed $y\in\Omega$ the function $x\mapsto
{\cal R}_j(x,y)$ is a solution of the boundary value problem
\begin{equation}\label{eq4.8}
\Delta_x^2 {\cal R}_j (x,y)= f_j(x,y) \,\,{\mbox{in}}\,\,\Omega,
\quad f_j(\cdot,y)\in C_0^{\infty}(\Omega),\quad {\cal
R}_j(\cdot,y)\in \ring W_2^2(\Omega),
\end{equation}

\noindent   so that
\begin{equation}\label{eq4.9}
\left\|\nabla_x^2{\cal
R}_j(\cdot,y)\right\|_{L^2(\Omega)}=\left\|{\cal
R}_j(\cdot,y)\right\|_{W_2^2(\Omega)}\leq
C\|f_j(\cdot,y)\|_{W^2_{-2}(\Omega)},\qquad j=1,2,3.
\end{equation}

\noindent Here $W^2_{-2}(\Omega)$ stands for the Banach space dual
of $\ring W_2^2(\Omega)$, i.e.
\begin{equation}\label{eq4.10}
\|f_j(\cdot,y)\|_{W^{2}_{-2}(\Omega)}=\sup_{v:\,\|v\|_{\ring
W_2^2(\Omega)}=1}\int_{\Omega}f_j(x,y)v(x)\,dx.
\end{equation}

 Recall that by Hardy's inequality
\begin{equation}\label{eq4.11}
\left\|\frac{v}{|\cdot-\,Q|^2}\right\|_{L^2(\Omega)}\leq C
\left\|\nabla^2 v\right\|_{L^2(\Omega)}\quad\mbox{for every}\quad
v\in\ring W_2^2(\Omega),\quad Q\in\po.
\end{equation}

\noindent Then for some $y_0\in\po$ such that $|y-y_0|=d(y)$
\begin{eqnarray}\label{eq4.12}\nonumber
&&\int_{\Omega}f_j(x,y)v(x)\,dx\leq C
\left\|\frac{v}{|\cdot-y_0|^2}\right\|_{L^2(\Omega)}\left\|f_j(\cdot,y)
|\cdot-y_0|^2\right\|_{L^2(\Omega)} \\[4pt]
&&\qquad \leq C d(y)^2 \left\|\nabla^2 v\right\|_{L^2(\Omega)}
\|f_j(\cdot,y) \|_{L^2(C_{d(y)/4, d(y)/2}(y))},
\end{eqnarray}

\noindent and therefore, by (\ref{eq4.7})
\begin{equation}\label{eq4.13}
\left\|\nabla_x^2{\cal R}(\cdot,y)\right\|_{L^2(\Omega)}\leq C
d(y)^{-1/2}.
\end{equation}

Turning to (\ref{eq4.2}), let us  first consider the case
$|x-y|\geq N d(y)$ for some large $N$ to be specified later. As
before, we denote by $y_0$ some point on the boundary such that
$|y-y_0|=d(y)$. Then by (\ref{eq4.7}) the function $x\mapsto {\cal
R}(x,y)$ is biharmonic in $\Omega\setminus B_{3d(y)/2}(y_0)$.
Hence, by Proposition~\ref{p3.4} with $r=6d(y)$
\begin{equation}\label{eq4.14}
|\nabla_x {\cal R}(x,y)|^2\leq
\frac{C}{|x-y_0|^2\,d(y)^3}\,\int_{C_{3d(y)/2,24d(y)}(y_0)}|{\cal
R}(z,y)|^2\,dz,
\end{equation}

\noindent provided $|x-y|\geq 4r+d(y)$, i.e $N\geq 25$. The
right-hand side of (\ref{eq4.14}) is bounded by
\begin{equation}\label{eq4.15}
\frac{Cd(y)}{|x-y_0|^2}\,\int_{C_{3d(y)/2,24d(y)}(y_0)}\frac{|{\cal
R}(z,y)|^2}{|z-y_0|^4}\,dz\leq \frac{C
\,d(y)}{|x-y_0|^2}\,\int_{\Omega}|\nabla_z^2{\cal R}(z,y)|^2\,dz
\leq \frac{C}{|x-y|^2},
\end{equation}

\noindent by Hardy's inequality and (\ref{eq4.13}).

 Now one can
directly check that
\begin{equation}\label{eq4.16}
|\nabla_x\nabla_y \Gamma(x,y)|\leq \frac{C}{|x-y|}\qquad\mbox{for
all}\qquad x,y\in\Omega,
\end{equation}

\noindent and combine it with  (\ref{eq4.14})--(\ref{eq4.15}) to
deduce that
\begin{equation}\label{eq4.17}
\bigg|\nabla_x\nabla_y (G(x,y)-\Gamma(x-y))\bigg|\leq
\frac{C}{|x-y|}\qquad\mbox{whenever}\qquad |x-y|\geq N d(y).
\end{equation}

\noindent  We claim that this settles the case
\begin{equation}\label{eq4.18}
|x-y|\geq N\min\{d(y),d(x)\}.
\end{equation}

\noindent Indeed, if $d(y)\leq d(x)$, (\ref{eq4.17}) gives the
desired result  and if $d(y)\geq d(x)$ and $|x-y|\geq N d(x)$, we
employ the version of (\ref{eq4.17}) with $d(x)$ in place of $d(y)$
which follows from the symmetry of Green's function and the
fundamental solution in $x$ and $y$ variables.

Next, assume that $|x-y|\leq N^{-1}d(y)$. For such $x$ we have
$\eta\bigl(\frac{x-y}{d(y)}\bigr)= 1$ and therefore
\begin{equation}\label{eq4.19}
\frac{\partial}{\partial y_j}\, (G(x,y)-\Gamma(x-y))={\cal
R}_j(x,y).
\end{equation}

\noindent By the interior estimates for solutions of elliptic
equations
\begin{equation}\label{eq4.20}
|\nabla_x{\cal R}(x,y)|^2\leq
\frac{C}{d(y)^5}\int_{B_{d(y)/8}(x)}|{\cal R}(z,y)|^2\,dz,
\end{equation}

\noindent since the function ${\cal R}$ is biharmonic in
$B_{d(y)/8}(x)\subset B_{d(y)/4}(y)$. Now we  bound the expression
above by
\begin{equation}\label{eq4.21}
\frac{C}{d(y)}\int_{B_{d(y)/4}(y)}\frac{|{\cal
R}(z,y)|^2}{|z-y_0|^4}\,dz \leq \frac{C}{d(y)}\left\|\nabla_x^2
{\cal R}(\cdot,y)\right\|_{L^2(\Omega)}^2\leq \frac{C}{d(y)^2}.
\end{equation}

\noindent When $|x-y|\leq N^{-1}d(y)$, we have
\begin{equation}\label{eq4.22}
(N-1)\, d(y)\leq N d(x)\leq (N+1)\, d(y),
\end{equation}

\noindent i.e. $d(y)\approx d(x)$, and therefore
(\ref{eq4.20})--(\ref{eq4.21}) give the desired result. By
symmetry, one can handle the case $|x-y|\leq N^{-1}d(x)$ and hence
all $x,y\in\Omega$ such that
\begin{equation}\label{eq4.23}
|x-y|\leq N^{-1}\,\max\{d(x),\,d(y)\}.
\end{equation}

Finally, it remains to consider the situation when
\begin{equation}\label{eq4.24}
|x-y|\approx d(x)\approx d(y),
\end{equation}

\noindent or more precisely, when
\begin{equation}\label{eq4.25}
N^{-1}\,d(x)\leq |x-y|\leq Nd(x)\quad \mbox{and}\quad
N^{-1}\,d(y)\leq |x-y|\leq Nd(y).
\end{equation}

In this case we use the biharmonicity of $x\mapsto G(x,y)$ in
$B_{d(x)/(2N)}(x)$. By the interior estimates, with $x_0\in\po$
such that $|x-x_0|=d(x)$, we have
\begin{eqnarray}\label{eq4.26}\nonumber
|\nabla_x\nabla_y G(x,y)|^2 &\leq &
\frac{C}{d(x)^5}\,\int_{B_{d(x)/(2N)}(x)}|\nabla_y
G(z,y)|^2\,dz\\[4pt]
&\leq&  \frac{C}{d(x)^5}\,\int_{B_{d(x)/(2N)}(x)}|\nabla_y
\Gamma(z-y)|^2\,dz+
\frac{C}{d(x)}\,\int_{B_{2d(x)}(x_0)}\frac{|{\cal
R}(z,y)|^2}{|z-x_0|^4}\,dz\nonumber \\[4pt]
&\leq&  \frac{C}{d(x)^5}\,\int_{B_{d(x)/(2N)}(x)}|\nabla_y
\Gamma(z-y)|^2\,dz+ \frac{C}{d(x)}\,\int_{\Omega}|\nabla_z^2 {\cal
R}(z,y)|^2\,dz\nonumber \\[4pt]
&\leq & \frac{C}{d(x)^{2}}+\frac{C}{d(x)d(y)},
\end{eqnarray}

\noindent invoking Hardy's inequality and (\ref{eq4.13}). In view of
(\ref{eq4.24}) this finishes the argument. \ep

\vskip 0.08in

\noindent {\it Proof of Theorem~\ref{t1.1}}. The estimate
(\ref{eq1.3}) follows directly from (\ref{eq4.2}). The second
inequality in (\ref{eq1.3.1}) can be proved closely following the
above argument, via an analogue of (\ref{eq4.2}). The first
inequality in (\ref{eq1.3.1}) is based on the second one and the
symmetry of Green's function. \ep


\vskip 0.08in  Green's function estimates proved in this section
allow to investigate the solutions of the Dirichlet problem
(\ref{eq1.1}) for a wide class of data. For example, consider the
boundary value problem
\begin{equation}\label{eq4.27}
\Delta^2 u={\rm div}\,  f+h, \quad u\in \ring W_2^2(\Omega),
\end{equation}

\noindent where $f=(f_1,f_2,f_3)$ is some vector valued function and
$h\in L^1(\Omega)$. Then the solution satisfies the estimate
\begin{equation}\label{eq4.28}
|\nabla u(x)|\leq C \int_\Omega \frac{|f(y)|}{|x-y|}\,dy+C
\int_\Omega |h(y)|\,dy,\qquad x\in\Omega.
\end{equation}

Indeed, the integral representation formula
\begin{equation}\label{eq4.29}
u(x)=\int_{\Omega} G(x,y)\,\Bigl({\rm div}\, f(y)+h(y)\Bigr)
\,dy,\qquad x\in\Omega,
\end{equation}

\noindent follows directly from the definition of Green's function.
It implies that
\begin{eqnarray}\label{eq4.30}\nonumber
\nabla u(x)&=& \nabla_x \int_{\Omega}  G(x,y)\,\Bigl({\rm div}\, f(y)+h(y)\Bigr)\,dy\\[4pt]
&=& -\int_{\Omega} \nabla_x (\nabla_y G(x,y)\cdot
f(y))\,dy+\int_{\Omega} \nabla_x  G(x,y) h(y)\,dy,
\end{eqnarray}

\noindent and Theorem~\ref{t1.1} leads to (\ref{eq4.28}).

One can further observe that by the mapping properties of the
Riesz potential the estimate (\ref{eq4.28}) entails that
\begin{equation}\label{eq4.31}
\|\nabla u\|_{L^\infty(\Omega)}\leq
C\|f\|_{L^{3/2,1}(\Omega)}+C\|h\|_{L^{1}(\Omega)},
\end{equation}

\noindent where $L^{3/2,1}(\Omega)$ is a Lorentz space.
Consequently,
\begin{equation}\label{eq4.32}
\|\nabla u\|_{L^\infty(\Omega)}\leq C\|
f\,\|_{L^{p}(\Omega)}+C\|h\|_{L^{1}(\Omega)}, \quad p>3/2,
\end{equation}

\noindent whenever $f\in L^p(\Omega)$ for some $p>3/2$.

\section{The capacity ${\rm
Cap}_P$} \setcounter{equation}{0}

This section is devoted to basic properties of the capacity ${\rm
Cap}_P$. A part of the results presented here and in \S{9} have been
obtained in \cite{MT}. For the convenience of the reader we present
a self-contained discussion.

To begin, we introduce a capacity of a compactum $K$ relative to
some open set $\Omega\subset\RR^3\setminus \{O\}$, $K\subset
\Omega$. To this end, recall that $\Pi$ is the space of functions
(\ref{eq1.7}) equipped with some norm. For example, we can take
\begin{equation}\label{eq5.1}
\|P\|_{\Pi}=\sqrt{b_0^2+b_1^2+b_2^2+b_3^2},
\end{equation}

\noindent and $\Pi_1:=\{P\in\Pi:\,\|P\|_{\Pi}=1\}$. A different
norm in the space $\Pi$ would yield an equivalent relative
capacity.

Now fix some $P\in\Pi_1$. Then
\begin{eqnarray}\label{eq5.2}
{\rm Cap}_P\,(K,\Omega):=\inf\Bigg\{\int_\Omega (\Delta
u(x))^2\,dx:\,\,u\in \ring W_2^2(\Omega),\,\,u=P\mbox{ in a
neighborhood of }K\Bigg\},
\end{eqnarray}

\noindent and
\begin{equation}\label{eq5.3}
{\rm Cap}\,(K,\Omega):=\inf_{P\in \Pi_1} {\rm Cap}_P\,(K,\Omega).
\end{equation}

\noindent Observe that in the introduction, for the sake of
brevity, we dropped the reference to $\Omega$. There we had
$\Omega=\RR^3\setminus\{0\}$.

It follows directly from the definition that the capacity ${\rm
Cap}_P$ is monotone in the sense that for every $P\in\Pi_1$
\begin{eqnarray}\label{eq5.4}
K_1\subseteq K_2\subset \Omega \quad & \Longrightarrow & \quad{\rm
Cap}_P\,(K_1,\Omega)\leq {\rm Cap}_P\,(K_2,\Omega),\\[4pt]\label{eq5.5}
K\subset \Omega_1\subseteq \Omega_2 \quad & \Longrightarrow &
\quad{\rm Cap}_P\,(K,\Omega_1)\geq {\rm Cap}_P\,(K,\Omega_2),
\end{eqnarray}

\noindent and analogous statements hold for ${\rm Cap}$ in place
of ${\rm Cap}_P$.

 We shall be concerned mostly with the case when a compactum is contained
in some annulus centered at the origin for the reasons that will
become apparent in the sequel. In such a case, it will be convenient
to work with an equivalent definition of capacity by means of the
form
\begin{equation}\label{eq5.6}
\Psi[u;\Omega]=\int_{\widetilde\varkappa(\Omega)} \left(
(\partial^2_r v)^2+2 r^{-2}(\partial_r v)^2+2 r^{-2}|\partial_r
\nabla_\omega v|^2+r^{-4}(\delta_\omega v)^2+2r^{-4}v\delta_\omega
v\right)\,r^2\,d\omega dr,
\end{equation}

\noindent where $(r,\omega)$ are the spherical coordinates in the
three dimensional space, $\tilde\varkappa$ is the mapping
\begin{equation}\label{eq5.7}
\RR^3 \ni x\,\stackrel{\widetilde\varkappa}{\longrightarrow}
\,(r,\omega)\in [0,\infty)\times S^{2},
\end{equation}

\noindent and $v=u\circ \widetilde\varkappa^{-1}$.

\begin{lemma}\label{l5.1} For every $r,R$ such that $0<r<R<\infty$ and every
function $u\in W_2^2(C_{r,R})$
\begin{eqnarray}\label{eq5.8}\nonumber
&&\hskip -1cm \Psi[u;C_{r,R}] =\int_{C_{r,R}} \Biggl[(\Delta u)^2\\[4pt]
 &&\hskip -1cm \qquad
 -\frac{2}{|x|^4}
\Bigl(x_i\frac{\partial}{\partial{x_i}}-1\Bigr) \left(\Bigl(x_j
\frac{\partial u}{\partial{x_j}}+u\Bigr) \Bigl( |x|^2 \Delta
u-x_i\frac{\partial}{\partial{x_i}}\Bigl(x_j\frac{\partial
u}{\partial{x_j}}\Bigr)-u\Bigr)+u^2\Bigr)\right) \Biggr]\,dx,
\end{eqnarray}

\noindent where, as customary, we imply summation on repeated
indices. Furthermore, for every open set $\Omega$ in
$\RR^3\setminus\{0\}$ and every $u\in \ring W_2^2(\Omega)$
\begin{equation}\label{eq5.9}\nonumber
\Psi[u;\Omega] =\int_{\Omega} (\Delta u(x))^2\,dx.
\end{equation}

\end{lemma}

The formulas (\ref{eq5.8})--(\ref{eq5.9}) can be checked directly
using the representation of the Laplacian in spherical coordinates
\begin{equation}\label{eq5.10}
\Delta
u=r^{-2}\left(\partial_r(r^2\partial_r)+\delta_\omega\right).
\end{equation}

They give rise to an alternative definition of the biharmonic
capacity. Indeed, if $K$ is a compact subset of $\Omega\subset
\RR^3\setminus\{0\}$, then for every $P\in\Pi_1$
\begin{eqnarray}\label{eq5.11}
{\rm Cap}_P\,(K,\Omega)=\inf\{\Psi[u; \Omega]:\,\,u\in \ring
W_2^2(\Omega),\,\,u=P\mbox{ in a neighborhood of }K\}
\end{eqnarray}

\noindent and an analogous equality holds for ${\rm Cap}$ in place
of ${\rm Cap}_P$.

\begin{lemma}\label{l5.2}
Suppose $K$ is a compactum in $\overline{C_{s,as}}$ for some
$s>0$, $a>1$. Then for every $P\in\Pi_1$
\begin{equation}\label{eq5.12}
{\rm Cap}_P(K, \RR^3\setminus\{O\})\approx  {\rm Cap}_P(K,
C_{s/2,2as})\qquad\mbox{and}\qquad {\rm Cap}_P(K, C_{s/2,2as})\leq
Cs^{-1},
\end{equation}

\noindent with the constants independent of $s$.
\end{lemma}

\bp The inequality
\begin{equation}\label{eq5.13} {\rm Cap}_P(K,
\RR^3\setminus\{O\})\leq {\rm Cap}_P(K,
C_{s/2,2as})
\end{equation}

\noindent is a consequence of the monotonicity property
(\ref{eq5.5}). As for the opposite inequality, we take $u\in \ring
W_2^2(\RR^3\setminus\{O\})$ such that $u=P$ in a neighborhood of
$K$ and
\begin{equation}\label{eq5.14}
{\rm Cap}_P(K, \RR^3\setminus\{O\})+\eps>\int_{\RR^3}|\Delta
u(x)|^2\,dx=\Psi[u;\RR^3\setminus\{O\}].
\end{equation}

\noindent Consider now the cut-off function
\begin{equation}\label{eq5.15}
\zeta \in C_0^\infty(1/2,2a), \,\,
\zeta=1\,\,\mbox{on}\,\,[3/4,3a/2],
\end{equation}

\noindent and let $w(x):=\zeta(|x|/s)u(x)$, $x\in \RR^3$. Then
\begin{equation}\label{eq5.16}
w\in\ring W_2^2(C_{s/2,2as})\quad\mbox{and}\quad w=P\,\,\mbox{in a
neighborhood of }\,\, K.
\end{equation}

\noindent Hence,
\begin{equation}\label{eq5.17}
{\rm Cap}_P(K, C_{s/2,2as})\leq \Psi[w; C_{s/2,2as}]
\end{equation}

\noindent and
\begin{eqnarray}\label{eq5.18}\nonumber
&&\Psi[w, C_{s/2,2as}]= \int_{s/2}^{2as}\int_{S^2} \Bigg(
(\partial^2_r (\zeta(r/s)v))^2+2 r^{-2}(\partial_r
(\zeta(r/s) v))^2\\[4pt]
&& \quad +2 r^{-2}|\partial_r (\zeta(r/s) \nabla_\omega v)|^2+
r^{-4}\zeta^2(r/s) (\delta_\omega v)^2+2r^{-4}\zeta^2(r/s)
v\delta_\omega v\Bigg)\,r^2\,d\omega dr\nonumber\\[4pt]
&& \quad \leq C\Psi[v, C_{s/2,2as}],
\end{eqnarray}

\noindent using the properties of $\zeta$ and the one dimensional
Hardy's inequality in the $r$ variable. This finishes the proof of
the first assertion in (\ref{eq5.12}).

As for the second one, observe first that if $v(x)=u(sx)$,
$x\in\RR^3$, the functions $u$ and $v$ belong to $\ring
W_2^2(\RR^3\setminus\{O\})$ simultaneously, and $u=P$ in a
neighborhood of $K$ if and only if $v=P$ in a neighborhood of
$s^{-1}K:=\{x\in\RR^3:\,sx\in K\}$. Also,
\begin{equation}\label{eq5.19}
\int_{\RR^3}|\Delta v(x)|^2\,dx= \int_{\RR^3}|\Delta_x
u(sx)|^2\,dx=s \int_{\RR^3}|\Delta_y u(y)|^2\,dy,
\end{equation}

\noindent so that
\begin{equation}\label{eq5.20}
s {\rm Cap}_P(K, \RR^3\setminus\{O\})= {\rm Cap}_P(s^{-1}K,
\RR^3\setminus\{O\}).
\end{equation}

\noindent However, $s^{-1} K\subset \overline{C_{1,a}}$ and
therefore by (\ref{eq5.12}) the right-hand side of (\ref{eq5.20})
is controlled by ${\rm Cap}_P(\overline{C_{1,a}},
\RR^3\setminus\{O\})$, uniformly in $s$. \ep

\begin{lemma}\label{l5.3}  Assume that for some $s>0$, $a>1$ the function $u\in
L^2(C_{s,as})$ is such that $\Psi[u; C_{s,as}]<\infty$. Then there
exists ${\cal P}={\cal P}(u,s,a)\in \Pi$ with the property
\begin{equation}\label{eq5.21}
\|u-{\cal P}\|_{L^2(C_{s,as})}^2\leq C s^4 \Psi[u; C_{s,as}].
\end{equation}
\end{lemma}

\bp Let us start with the expansion of $u$ by means of spherical
harmonics:
\begin{equation}\label{eq5.22}
u=\sum_{l=0}^{\infty}\sum_{m=-l}^{l}u_l^m(r) Y_l^m(\omega),
\end{equation}

\noindent where $Y_l^m$ are the spherical harmonics of degree
$l\in\NN$ and order $m\in\ZZ$.
By Poincar{\'e}'s inequality, for $l=0,1$, and the corresponding
$m$ there exist constants $\sigma_l^m$ (depending on $u$) such
that
\begin{equation}\label{eq5.23}
\int_{s}^{as}|u_l^m(r)- \sigma_l^m|^2\,dr \leq C s^2
\int_{s}^{as}|\partial_r u_l^m(r)|^2\,dr.
\end{equation}

\noindent Let
\begin{equation}\label{eq5.24}
{\cal P}(x):=\sigma_0^0
+{\sigma_1^1}\,\frac{x_1}{|x|}+{\sigma_1^{-1}}\,\frac{x_2}{|x|}+{\sigma_1^0}\,\frac{x_3}{|x|},
\qquad x\in\RR^3\setminus \{O\}.
\end{equation}

\noindent Then (\ref{eq5.23}) yields (\ref{eq5.21}). \ep



\begin{proposition}\label{l5.4} (\cite{MT}) Suppose $s>0$, $a\geq 2$ and $K$ is a
compact subset of $\overline{C_{s,as}}$. Then for every $u\in
L^2(C_{s,as})$ such that $\Psi[u; C_{s,as}]<\infty$ and $u=0$ in a
neighborhood of $K$
\begin{equation}\label{eq5.25}
\frac{1}{s^3}\int_{C_{s,as}}|u(x)|^2\,dx\leq \frac{C}{{\rm
Cap}\,(K,\RR^3\setminus \{O\})} \,\,\Psi[u; C_{s,as}].
\end{equation}
\end{proposition}

\bp For the purposes of this argument let us take
$\|P\|_{\Pi}:=\|P\|_{L^2(C_{1,a})}$ and  let
$\Pi_1:=\{P\in\Pi:\,\|P\|_{\Pi}=1\}$ with such a norm.  This is an
equivalent norm in the space $\Pi$ and hence it yields the
capacity equivalent to the one defined in
(\ref{eq5.1})--(\ref{eq5.2}). We claim that for every $P\in\Pi_1$
\begin{equation}\label{eq5.26}
{\rm Cap}_P(K, C_{s/2,2as})\leq Cs^{-4}\|P-u\|_{L^2(C_{s,as})}^2+C
\Psi[u; C_{s,as}].
\end{equation}

To prove this, let us denote by $V^2_2(C_{s,as})$ a collection of
functions on $C_{s,as}$ such that
\begin{equation}\label{eq5.27}
\|u\|_{V^2_2(C_{s,as})}:=\left(\frac{1}{s^4}\int_{C_{s,as}}|u(x)|^2\,dx+
\Psi[u; C_{s,as}]\right)^{1/2},
\end{equation}

\noindent is finite. One can construct an extension operator
\begin{equation}\label{eq5.28}
{\rm Ex}: V^2_2(C_{s,as})\to V^2_2(C_{s/2,2as})
\end{equation}

\noindent with the operator norm independent of $s$ satisfying the
properties
\begin{equation}\label{eq5.29}
{\rm Ex}\, u=u\mbox{ in }C_{s,as},\qquad {\rm Ex}\, P=P\mbox{ for
every }P\in\Pi_1,
\end{equation}

\noindent and such that if $u=0$ in some neighborhood of $K$
intersected with $\overline{C_{s,as}}$ then ${\rm Ex}\, u$
vanishes in a neighborhood of $K$ contained in $C_{s/2,2as}$. For
example, one can start with the corresponding one-dimensional
extension operator and then use the expansion (\ref{eq5.22}) to
define ${\rm Ex}$.

Having this at hand, we define $w(x):=\zeta(|x|/s)(P(x)-{\rm Ex}\,
u(x))$, $x\in C_{s/2,2as}$, where $\zeta$ is a function introduced
in  (\ref{eq5.15}). Then $w$ satisfies (\ref{eq5.16}) and
therefore ${\rm Cap}_P(K, C_{s/2,2as})$ is controlled by
\begin{eqnarray}\label{eq5.30}\nonumber
&&\Psi[w; C_{s/2,2as}]\leq \Psi[P-{\rm Ex}\, u;
C_{s/2,2as}]=\Psi[{\rm Ex}\,(P-u); C_{s/2,2as}]\nonumber\\[4pt]
&& \qquad \leq Cs^{-4}\|P-u\|_{L^2(C_{s,as})}^2+C \Psi[P-u;
C_{s,as}],
\end{eqnarray}

\noindent where the first inequality is proved analogously to
(\ref{eq5.18}) and the second one follows from the mapping
properties of ${\rm Ex}$. Using that $\delta_\omega
\omega_i=-2\omega_i$, $i=1,2,3,$ one can directly check that
\begin{equation}\label{eq5.31}
\Psi[P-u; C_{s,as}]=\Psi[u; C_{s,as}],
\end{equation}

\noindent and obtain  (\ref{eq5.26}).

The next step is to pass from (\ref{eq5.26}) to (\ref{eq5.25}).
Without loss of generality we may assume that
$\|u\|_{L^2(C_{s,as})}=s^{3/2}$. Then the desired result reads as
\begin{equation}\label{eq5.32}
\inf_{P\in\Pi_1} {\rm Cap}_P(K, C_{s/2,2as})\leq \Psi[u;
C_{s,as}].
\end{equation}

Let ${\cal P}={\cal P}(u,s,a)$ be a function in $\Pi$ satisfying
(\ref{eq5.21}), and denote by $C_0$ the constant $C$ in
(\ref{eq5.21}). First of all, the case
\begin{equation}\label{eq5.33}
\Psi[u; C_{s,as}]\geq 1/(4C_0s)
\end{equation}

\noindent is trivial, since Lemma~\ref{l5.2} guarantees that the
right-hand side of (\ref{eq5.33}) is bounded from below by the
capacity of $K$, modulo a multiplicative constant.

On the other hand,
\begin{equation}\label{eq5.34}
\Psi[u; C_{s,as}]\leq 1/(4C_0s)\qquad\Longrightarrow \qquad 2
\|u-{\cal P}\|_{L^2(C_{s,as})}\leq s^{3/2}=\|u\|_{L^2(C_{s,as})},
\end{equation}

\noindent by (\ref{eq5.21}) and the normalization of $u$. This, in
turn, implies that
\begin{equation}\label{eq5.35}
\frac{s^{3/2}}{2}\leq \|{\cal P}\|_{L^2(C_{s,as})}\leq
\frac{3s^{3/2}}{2}.
\end{equation}

\noindent Finally, we choose
\begin{equation}\label{eq5.36}
P:=\frac{{\cal P}}{\|{\cal
P}\|_{L^2(C_{1,a})}}=s^{3/2}\,\frac{{\cal P}}{\|{\cal
P}\|_{L^2(C_{s,as})}}.
\end{equation}

\noindent Then
\begin{eqnarray}\label{eq5.36.1}
&&\|P-{\cal P}\|_{L^2(C_{s,as})} = \left| s^{3/2} - \|{\cal
P}\|_{L^2(C_{s,as})}
\right|\nonumber\\[4pt]&&\qquad =\left| \|u\|_{L^2(C_{s,as})} - \|{\cal P}\|_{L^2(C_{s,as})}
\right| \leq \| u - {\cal P} \|_{L^2(C_{s,as})}.
\end{eqnarray}

\noindent Hence,
\begin{equation}\label{eq5.36.2}
||u-P|| \leq ||u-{\cal P}|| + ||{\cal P}-P|| \leq 2 ||u-{\cal
P}||,
\end{equation}

\noindent so that
\begin{equation}\label{eq5.37}
\|u-P\|_{L^2(C_{s,as})}^2\leq 16 \|u-{\cal P}\|_{L^2(C_{s,as})}^2
\leq 16C_0s^4 \Psi[u; C_{s,as}],
\end{equation}

\noindent by (\ref{eq5.21}). Combining (\ref{eq5.37}) with
(\ref{eq5.26}), we complete the argument. \ep

\section{1-regularity of a boundary point}
\setcounter{equation}{0}

Let $\Omega$ be a domain in $\RR^3$ and consider the boundary
value problem
\begin{equation}\label{eq6.1}
\Delta^2 u=f \,\,{\mbox{in}}\,\,\Omega, \quad f\in
C_0^{\infty}(\Omega),\quad u\in \ring W_2^2(\Omega).
\end{equation}

\noindent We say that the point $Q\in\po$ is {\it 1-regular} (with
respect to $\Omega$) if for every $f\in C_0^{\infty}(\Omega)$ the
gradient of the solution to (\ref{eq6.1}) is continuous, i.e.
\begin{equation}\label{eq6.2}
\nabla u(x)\to 0\mbox{ as } x\to Q,\,x\in\Omega.
\end{equation}

\noindent Otherwise, $Q\in\po$ is called 1-irregular.

Observe that in the case $Q=O$ this definition coincides with the
one given in the introduction.

In this section we would like to show that 1-regularity is a local
property. In particular, while most of the statements in Sections
1--5 were confined to the case of a bounded domain, the
proposition below will allow us to study 1-regularity with respect
to any open set in $\RR^3$.

\begin{proposition}\label{p6.1}
Let $\Omega$ be a bounded domain in $\RR^3$ and the point
$Q\in\po$ be 1-regular with respect to $\Omega$. If $\Omega'$ is
another domain with the property that $B_r(Q)\cap\Omega =
B_r(Q)\cap\Omega'$ for some $r>0$ then $Q$ is 1-regular with
respect to $\Omega'$.
\end{proposition}

The proof of the proposition rests on the ideas from \cite{M2}. It
starts with the following result.

\begin{lemma}\label{l6.2}
Let $\Omega$ be a bounded domain in $\RR^3$ and the point
$Q\in\po$ be 1-regular with respect to $\Omega$. Then
\begin{equation}\label{eq6.3}
\nabla u(x)\to 0\mbox{ as } x\to Q,\,x\in\Omega,
\end{equation}

\noindent for every $u\in \ring W_2^2(\Omega)$ satisfying
\begin{equation}\label{eq6.4}
\Delta^2 u=\sum_{\alpha:\,|\alpha|\leq 2} \partial^\alpha f_\alpha
\,\,{\mbox{in}}\,\,\Omega, \quad f_\alpha\in L^2(\Omega)\cap
C^\infty(\Omega),\quad f_\alpha=0\mbox{ in a neighborhood of } Q.
\end{equation}
\end{lemma}

\bp Take some $\eta\in C_0^\infty(\Omega)$ and let $v$ be the
solution of the Dirichlet problem
\begin{equation}\label{eq6.5}
\Delta^2 v=\sum_{\alpha:\,|\alpha|\leq 2} \partial^\alpha (\eta
f_\alpha) \,\,{\mbox{in}}\,\,\Omega, \quad v\in \ring
W_2^2(\Omega),
\end{equation}

\noindent and $w:=u-v\in W_2^2(\Omega)$. Since the point $Q$ is
1-regular, the function $v$ automatically satisfies (\ref{eq6.3})
and it remains to consider $w$.

Since $f_\alpha=0$ in a neighborhood of $Q$, the function $w$ is
biharmonic in some neighborhood of $Q$ and, therefore, for some
$R>0$ depending on the ${\rm supp}\,f_\alpha$, we have
\begin{equation}\label{eq6.6}
|\nabla w(x)|\leq  \frac{C}{d(x)^3} \int_{B_{ d(x)/2}(x)} |\nabla
w(y)|^2\,dy \leq \frac{C}{R^{5}}\int_{C_{R/4,4R}(Q)}
|w(y)|^2\,dy,\quad \forall\,\,x\in B_{R/4}(Q),
\end{equation}

\noindent analogously to (\ref{eq3.17})--(\ref{eq3.18}). On the
other hand, according to Lemma~\ref{l2.3} the last expression in
(\ref{eq6.6}) does not exceed
\begin{eqnarray}\label{eq6.7}\nonumber
&& \hskip -1cm C\,\sup_{\xi\in C_{R/4,4R}(Q)\cap\Omega}
\int_{\RR^n}\Delta w(y)\,\Delta\Bigg(\frac{w(y)}{|x-Q|}\,g\bigg(\log
\frac{|\xi-Q|}{|y-Q|}\bigg)\Bigg)\,dy\\[4pt]
&& \hskip -1cm \leq C \,\sup_{\xi\in C_{R/4,4R}(Q)\cap\Omega}
\sum_{\alpha:\,|\alpha|\leq 2}
\int_{\RR^n}(1-\eta(y))f_\alpha(y)\,\,(-\partial_y)^\alpha
\Bigg(\frac{w(y)}{|y-Q|}\,g\bigg(\log
\frac{|\xi-Q|}{|y-Q|}\bigg)\Bigg)\,dy,
\end{eqnarray}

\noindent where $g$ is given by (\ref{eq2.8}). When $x$ approaches
$Q$, the support of $1-\eta$ can be chosen arbitrarily small.
Hence, the integral on the right-hand side of  (\ref{eq6.7})
shrinks. Then (\ref{eq6.6})--(\ref{eq6.7}) imply that $|\nabla
w(x)|\to 0$ when $x\to Q$. \ep

\vskip 0.08in

\noindent {\it Proof of Proposition~\ref{p6.1}.}\, Consider a
solution of the Dirichlet problem
\begin{equation}\label{eq6.8}
\Delta^2 u=f \,\,{\mbox{in}}\,\,\Omega', \quad f\in
C_0^{\infty}(\Omega'),\quad u\in \ring W_2^2(\Omega'),
\end{equation}

\noindent and take some cut-off function $\eta\in
C_0^\infty(B_{r}(Q))$ equal to 1 on $B_{r/2}(Q)$. Then $\eta u \in
\ring W_2^2(\Omega)$ and
\begin{equation}\label{eq6.9}
\Delta^2 (\eta u)=\eta f+[\Delta^2,\eta] u.
\end{equation}

Since $\eta f\in C_0^\infty(\Omega)$,
\begin{equation}\label{eq6.10}
[\Delta^2,\eta]:\ring W_2^2(\Omega)\longrightarrow (\ring
W_2^2(\Omega))^*=W_{-2}^2(\Omega)\quad {\mbox{and}}\quad {\rm
supp}\,([\Delta^2,\eta]u)\subset C_{r/2,r}(Q)\cap\Omega,
\end{equation}

\noindent one can write
\begin{equation}\label{eq6.11}
\Delta^2(\eta u)= \sum_{\alpha:\,|\alpha|\leq 2} \partial^\alpha
f_\alpha, \quad\mbox{for some}\quad f_\alpha\in L^2(\Omega)\cap
C^\infty(\Omega),
\end{equation}

\noindent with $f_\alpha=0$ in a neighborhood of $Q$ given by the
intersection of $B_{r/2}(Q)$ and the complement to ${\rm
supp}\,f$. Then, by Lemma~\ref{l6.2}, the gradient of $\eta u$
(and therefore, the gradient of $u$) vanishes as $x\to Q$. \ep

\section{Sufficient condition for 1-regularity}
\setcounter{equation}{0}


The following proposition provides the first part of
Theorem~\ref{t1.2}, i.e. sufficiency of condition (\ref{eq1.9})
for 1-regularity of a boundary point.
\begin{proposition}\label{p7.1}
Let $\Omega$ be a bounded domain in $\RR^3$,
$O\in\RR^3\setminus\Omega$, $R>0$ and
\begin{equation}\label{eq7.1}
\Delta^2 u=f \,\,{\mbox{in}}\,\,\Omega, \quad f\in
C_0^{\infty}(\Omega\setminus B_{4R}),\quad u\in \ring
W_2^2(\Omega).
\end{equation}

\noindent Fix some $a\geq 4$. Then for every $x\in
B_{R/a^4}\cap\Omega$
\begin{eqnarray}\label{eq7.2}\nonumber
|\nabla u(x)|^2+\frac{|u(x)|^2}{|x|^2}&\leq&
\frac{C}{R^{5}}\int_{C_{R,4R}\cap\Omega}
|u(y)|^2\,dy\,\,\\[4pt]&&
\times\exp \left( -c\sum_{j=2}^{l}(Ra^{-2j})\,{\rm
Cap}\,(\overline{C_{R \,a^{-2j},R\,a^{-2(j-1)}}}\setminus\Omega,\RR^3\setminus
\{O\})\right),
\end{eqnarray}
\noindent where $l\geq 2$, $l\in\NN$, is such that $|x|\leq a^{-2l}R$.

In particular, when $O$ is a boundary point of $\Omega$,
\begin{equation}\label{eq7.3}
\mbox{if } \sum_{j=1}^{\infty}a^{-j} {\rm
Cap}\,(\overline{C_{a^{-j},a^{-(j-1)}}}\setminus\Omega,\RR^3\setminus
\{O\})\,ds=+\infty\,\,\, \mbox{then $O$ is 1-regular}.
\end{equation}
Here $a$ is any real number greater than 1.
\end{proposition}

\bp Fix $s\leq {R}/{a^2}$ and let us introduce some extra
notation. First,
\begin{equation}\label{eq7.3.1}
\gamma(s):={{\rm
Cap}\,(\overline{C_{s,a^2s}}\setminus\Omega,\RR^3\setminus
\{O\})}.
\end{equation}

\noindent Further, let $Q_\tau [u;\Omega]$, $\tau\in\RR$, be the
quadratic form
\begin{eqnarray}\label{eq7.3.2}
Q_\tau [u;\Omega]&:=&\int_{\varkappa(\Omega)}\Bigl[(\delta_\omega
v)^2g(t-\tau)+2(\partial_t\nabla_\omega v)^2g(t-\tau)+
(\partial_t^2v)^2g(t-\tau)\nonumber\\[4pt]
&&  -(\nabla_\omega v)^2 \Bigl(\partial_t^2
g(t-\tau)+\partial_tg(t-\tau)+2g(t-\tau)\Bigr) \nonumber\\[4pt]
&& -(\partial_t v)^2\Bigl(2\partial_t^2g(t-\tau)+3\partial_t
g(t-\tau)-g(t-\tau)\Bigr)\Bigr]\,d\omega dt,
\end{eqnarray}

\noindent where $v=e^t(u\circ \varkappa^{-1})$, $g$ is defined by
(\ref{eq2.8}) and $\varkappa$ is the change of coordinates
(\ref{eq2.2}). Throughout this proof $\tau=\log s^{-1}$.

Now take  $\eta\in C_0^\infty(B_{2s})$ such that
\begin{equation}\label{eq7.3.3}
0\leq\eta\leq 1\,\, {\rm in}\,\,B_{2s}, \quad \eta=1\,\, {\rm
in}\,\,B_s\quad {\rm and} \quad |\nabla^k \eta|\leq C/|x|^k, \quad
k\leq 4.
\end{equation}

\noindent  Following the argument in (\ref{eq3.8})--(\ref{eq3.10})
and the discussion after (\ref{eq3.10}), and then passing to the
limit as $n\to\infty$, we have
\begin{eqnarray}\label{eq7.3.4}
Q_\tau [u;B_s]&\leq &Q_\tau [\eta u;\Omega]\leq \int_{\RR^3}\Delta
\bigg(\eta(x)u(x)\bigg)\, \Delta\bigg(\eta(x)u(x)|x|^{-1}g(\log
(s/|x|))\bigg) \,dx\nonumber\\[4pt]
&\leq & C\sum_{k=0}^2 \frac{1}{s^{5-2k}}\int_{C_{s,2s}} |\nabla^k
u(x)|^2\,dx\leq \frac{C}{s^{5}}\int_{C_{s,4s}} |u(x)|^2\,dx.
\end{eqnarray}

\noindent Denote
\begin{equation}\label{eq7.3.5}
\varphi(s):=\sup_{|x|\leq s}\left(|\nabla u(x)|^2 +
\frac{|u(x)|^2}{|x|^2}\right)+Q_\tau [u;B_s],\qquad \tau=\log
s^{-1},\qquad s\leq \frac{R}{a^2}.
\end{equation}

\noindent Then combining (\ref{eq7.3.4}) with Corollary~\ref{c3.3}
and Proposition~\ref{p3.2},
\begin{equation}\label{eq7.3.6}
\varphi(s)\leq \frac{C}{s^5} \int_{C_{s,16s}} |u(x)|^2\,dx\leq
\frac{C}{s^5} \int_{C_{s,a^2s}} |u(x)|^2\,dx.
\end{equation}

\noindent For $\gamma(s)>0$ the expression on the right-hand side
of (\ref{eq7.3.6}) does not exceed
\begin{equation}\label{eq7.5}
\frac{C}{s^3} \int_{C_{s,a^2s}} \frac{|u(x)|^2}{|x|^2}\,dx\leq
\frac{C}{\gamma(s)}\,\, \Psi\left[\frac{u}{|x|};C_{s,a^2s}\right]
\leq \frac{C}{s\gamma(s)}\,\,Q_\tau [u;C_{s,a^2s}],
\end{equation}

\noindent where we used Proposition~\ref{l5.4} for the first
inequality. The second one can be proved directly using that
$e^{-\tau}=s$ and calculations from the proof of Lemma~\ref{l2.3}.
All in all,
\begin{equation}\label{eq7.7}
\varphi(s)\leq \frac{C}{s\gamma(s)}\,\,Q_\tau
[u;C_{s,a^2s}]\leq\frac{C}{s\gamma(s)}\left(\varphi(a^2s)-\varphi(s)\right),
\end{equation}

\noindent which, in turn, implies that
\begin{equation}\label{eq7.8}
\varphi(s)\leq\frac{1}{1+C^{-1}\, s\gamma(s)}\,\varphi(a^2s)\leq
\exp\left(-cs\gamma(s)\right)\,\varphi(a^2s),
\end{equation}

\noindent since $s\gamma(s)$ is bounded by (\ref{eq5.12}). In
particular, employing (\ref{eq7.8}) for $s=a^{-2j}r$, $r\leq R$,
$j\in\NN$, one can conclude that
\begin{equation}\label{eq7.9}
\varphi(a^{-2l}r)\leq
\exp\left(-c\sum_{j=2}^la^{-2j}r\,\gamma(a^{-2j}r)\right)
\,\varphi(ra^{-2}),
\end{equation}

\noindent for all $l\geq 2$, $l\in\NN$.

Let us choose $l\geq 2$, $l\in\NN$ so that
\begin{equation}\label{eq7.10}
a^{-2l-2}R\leq |x|\leq a^{-2l}R.
\end{equation}
Using (\ref{eq7.9}), we deduce that for every $x\in B_{R/a^4}\cap\Omega$
and $l$ defined by (\ref{eq7.10})
\begin{equation}\label{eq7.12}
|\nabla u(x)|^2 + \frac{|u(x)|^2}{|x|^2} \leq
\varphi(a^{-2l}R)\leq
\exp\left(-c\sum_{j=2}^la^{-2j}R\,\gamma(a^{-2j}R)\right)
\,\varphi(R a^{-2}).
\end{equation}

 Finally, analogously to (\ref{eq7.3.4})--(\ref{eq7.3.6})
\begin{eqnarray}\label{eq7.14}\nonumber
\varphi(Ra^{-2})&\leq & \sup_{|x|\leq R/a^2}\left(|\nabla u(x)|^2 +
\frac{|u(x)|^2}{|x|^2}\right)+C\sum_{k=0}^2 \int_{C_{R/a^2,2R/a^2}}
\frac{|\nabla^ku(x)|^2}{|x|^{5-2k}}\,dy\\[4pt]&\leq&
\frac{C}{R^{5}}\int_{C_{R/a^2,16R/a^2}} |u(y)|^2\,dy\leq
\frac{C}{R^{5}}\int_{C_{R,4R}} |u(y)|^2\,dy,
\end{eqnarray}

\noindent using Proposition~\ref{p3.2} for the last inequality.
Therefore, we finish the proof
of (\ref{eq7.2}). 

Now let us turn to (\ref{eq7.3}). The estimate (\ref{eq7.2}) directly leads to the following conclusion. When $O$ is a boundary point of $\Omega\subset \RR^n$ 
\begin{equation}\label{eq7.3-aux}
\mbox{if } \sum_{j=1}^{\infty}(a^{-j}R) {\rm
Cap}\,(\overline{C_{a^{-j}R,a^{-(j-1)}R}}\setminus\Omega,\RR^3\setminus\{O\})=+\infty\,\,\, \mbox{then $O$ is 1-regular},
\end{equation}
where $a\geq 16$ (we swapped $a^2$ in (\ref{eq7.2}) for $a$ in (\ref{eq7.3-aux})). Next, the condition $a\geq 16$ can be substituted by any $a>1$ using monotonicity of capacity to shrink $\overline{C_{R \,a^{-j},R\,a^{-(j-1)}}}\setminus\Omega$ as necessary.  Finally, there exists $N\in\ZZ$ such that $R\approx a^{-N}$, so that the series in (\ref{eq7.3-aux}) can be rewritten as the series in (\ref{eq7.3}), with the summation over $j=N+1, N+2, ...$, but that again does not affect the question of convergence. Hence, we arrive at (\ref{eq7.3}).

\ep

Given the result of Proposition~\ref{p7.1}, we can derive the
estimates for biharmonic functions at infinity as well as those
for Green's function in terms of the capacity of the complement of
$\Omega$, in the spirit of (\ref{eq7.2}).

\begin{proposition}\label{p7.2} Let $\Omega$ be a bounded domain in $\RR^3$,
$O\in\RR^3\setminus\Omega$, $r>0$ and assume that
\begin{equation}\label{eq7.15}
\Delta^2 u=f \,\,{\mbox{in}}\,\,\Omega, \quad f\in
C_0^{\infty}(B_{r/4}\cap\Omega),\quad u\in \ring W_2^2(\Omega).
\end{equation}

\noindent Fix some $a\geq 4$. Then for any $x\in\Omega\setminus
B_{a^4r}$
\begin{eqnarray}\label{eq7.16}
\nonumber |\nabla u(x)|^2+\frac{|u(x)|^2}{|x|^2}&  \leq &
\frac{C}{|x|^2\,r^3}\int_{C_{\frac r4, r}\cap\Omega}
|u(y)|^2\,dy\,\,\\[4pt]&&
\times\exp \left( -c\sum_{j=2}^{l}(ra^{2j}) {\rm
Cap}\,(\overline{C_{ra^{2(j-1)},ra^{2j}}}\setminus\Omega,\RR^3\setminus\{O\})\right),
\end{eqnarray}
where $l\geq 2$, $l\in\NN$, is such that $|x|\geq a^{2l}r$.
\end{proposition}

\bp Recall the proof of Proposition~\ref{p3.4}. With the notation
(\ref{eq3.26}) the results (\ref{eq3.27})--(\ref{eq3.29}),
(\ref{eq3.32}) allow to apply  Proposition~\ref{p7.1} to $U_n$,
$R=1/r$, in order to write
\begin{eqnarray}
&&\hskip -1cm |\nabla u_n(x)|^2+\frac{|u_n(x)|^2}{|x|^2} \leq
C\,\frac{ \left|(\nabla
U_n)(x/|x|^2)\right|^2}{|x|^2}+\left|U_n(x/|x|^2)\right|^2\nonumber\\[4pt]\nonumber
&&\hskip -1cm \qquad  \leq C\,\frac{r^5}{|x|^2}\int_{C_{\frac 1r,
\frac 4r}} |U_n(z)|^2\,dz\,\,\\[4pt]
&&\times\exp \left( -c\sum_{j=2}^{l}(a^{-2j}/r) {\rm
Cap}\,(\overline{C_{a^{-2j}/r,a^{-2(j-1)}/r}}\setminus{\mathcal I}(\Omega_n),\RR^3\setminus
\{O\})\right).
\end{eqnarray}
Here $a\geq 4$, $l\geq 2$, $l\in\NN$, and $x$ is such that $|x|\geq a^{2l}r$.

We claim that
\begin{equation}\label{eq7.17}
 {\rm
Cap}\,(\overline{C_{\frac 1s,\frac {a^2}s}}\setminus{\cal I}(\Omega_n),\RR^3\setminus
\{O\})\approx s^2 \,{\rm
Cap}\,(\overline{C_{\frac{s}{a^2},s}}\setminus\Omega_n,\RR^3\setminus
\{O\}),
\end{equation}

\noindent where the implicit constants are independent of $s$.

Indeed,
\begin{equation}\label{eq7.18}
 {\rm
Cap}\,(\overline{C_{\frac 1s,\frac {a^2}s}}\setminus{\cal
I}(\Omega_n),\RR^3\setminus \{O\})\approx {\rm
Cap}\,(\overline{C_{\frac 1s,\frac {a^2}s}}\setminus{\cal
I}(\Omega_n),C_{\frac 1{2s},\frac {2a^2}{s}}),
\end{equation}

\noindent and for every $u\in \ring W_2^2(C_{\frac 1{2s},\frac
{2a^2}{s}})$ the function $y\mapsto |y|\, u(y/|y|^2)$ belongs to
$\ring W_2^2(C_{\frac s{2a^2},2s})$ by (\ref{eq3.29}) and therefore,
if $U(y):=u(y/|y|^2)$ then $U\in \ring W_2^2(C_{\frac s{2a^2},
2s})$. In addition, if $u=P$ in a neighborhood of
${\overline{C_{\frac 1s,\frac {a^2}s}}}\setminus{\cal I}(\Omega_n)$
then $U(y)=P(y/|y|^2)=P(y)$ for all $y$ in the corresponding
neighborhood of ${\overline{C_{\frac sa^2, s}}}\setminus\Omega_n$.
Finally, by (\ref{eq3.28})
\begin{equation}\label{eq7.19}
\int_{C_{\frac 1{2s},\frac {2a^2}{s}}} |\Delta u(x)|^2\,
dx=\int_{C_{\frac s{2a^2},2s}} |\Delta (|y|\, u(y/|y|^2))|^2\,
dy\approx s^2 \int_{C_{\frac s{2a^2},2s}} |\Delta U(y)|^2\, dy,
\end{equation}

\noindent since $u\in \ring W_2^2(C_{\frac 1{2s},\frac {2a^2}{s}})$.
This proves the ``$\geq$" inequality in (\ref{eq7.17}). The
opposite inequality reduces to the previous one taking $1/s$ in
place of $s$ and ${\cal I}(\Omega_n)$ in place of $\Omega_n$,
since ${\cal I}({\cal I}(\Omega_n))=\Omega_n$.

As a result, we have
\begin{eqnarray}\label{eq7.20}
&&\hskip -1cm |\nabla u_n(x)|^2+\frac{|u_n(x)|^2}{|x|^2}
\nonumber\\[4pt]
&&\hskip -1cm \qquad  \leq \frac{C}{|x|^2\,r^3}\int_{C_{\frac r4,
r}}|u_n(z)|^2\,dz \,\,\nonumber\\[4pt]
&& \times\exp \left( -c\sum_{j=2}^{l}(ra^{2j}) {\rm
Cap}\,(\overline{C_{ra^{2(j-1)},ra^{2j}}}\setminus\Omega_n,\RR^3\setminus\{O\})\right),
\end{eqnarray}

\noindent using the monotonicity property (\ref{eq5.4}). Now the
argument can be finished using the limiting procedure similar to
the one in Proposition~\ref{p3.4}. \ep

The following Proposition is a more precise version of the
estimate on Green's function we announced in the introduction
after Theorem~\ref{t1.2}.
\begin{proposition}\label{p7.3} Let $\Omega$ be a bounded domain in $\RR^3$,
$O\in\po$. Fix some $a\geq 4$ and let $c_a:=1/(32a^4)$. Then for
$x,y\in\Omega$
\begin{eqnarray}
\nonumber &&|\nabla_x\nabla_y G(x,y)|\\[4pt]\nonumber
&& \quad \leq  \left\{\begin{array}{l} \frac{C}{|x-y|}
\times\exp \left( -c\sum_{j=2}^{l_{yx}}(|y|a^{2j}) {\rm
Cap}\,(\overline{C_{32|y|a^{2(j-1)},32|y|a^{2j}}}\setminus\Omega,\RR^3\setminus\{O\})\right),\\[8pt]\nonumber
 \quad\quad\mbox{if  $|y|\leq c_a |x|$ and \,$l_{yx}\geq 2$, $l_{yx}\in\NN$, is such that $|x|\geq 32 a^{2l_{yx}}|y|$},\\[8pt]
\frac{C}{|x-y|} \times\exp \left( -c\sum_{j=2}^{l_{xy}}(|x|a^{2j}) {\rm
Cap}\,(\overline{C_{32|x|a^{2(j-1)},32|x|a^{2j}}}\setminus\Omega,\RR^3\setminus\{O\})\right),\\[8pt]\nonumber
 \quad\quad\mbox{if  $|x|\leq c_a |y|$ and \,$l_{xy}\geq 2$, $l\in\NN$, is such that $|y|\geq 32 a^{2l_{xy}}|x|$},\\[8pt]
\frac{C}{|x-y|}, \qquad \mbox{if}\quad c_a |y|\leq |x|\leq
c_a^{-1}|y|,
\end{array}
\right.
\end{eqnarray}

\noindent and
\begin{eqnarray}
\nonumber &&\max\Bigl\{|\nabla_x G(x,y)|, \,|\nabla_y G(x,y)|\Bigr\}\\[4pt]\nonumber
&& \quad \leq  \left\{\begin{array}{l}
C\,\exp \left( -c\sum_{j=2}^{l_{yx}}(|y|a^{2j}) {\rm
Cap}\,(\overline{C_{32|y|a^{2(j-1)},32|y|a^{2j}}}\setminus\Omega,\RR^3\setminus\{O\})\right),\\[8pt]\quad\mbox{if  $|y|\leq c_a |x|$ and \,$l_{yx}\geq 2$, $l_{yx}\in\NN$, is such that $|x|\geq 32a^{2l_{yx}}|y|$,}\\[8pt]
C\,\exp \left( -c\sum_{j=2}^{l_{xy}}(|x|a^{2j}) {\rm
Cap}\,(\overline{C_{32|x|a^{2(j-1)},32|x|a^{2j}}}\setminus\Omega,\RR^3\setminus\{O\})\right),\\[8pt]\nonumber
 \quad\quad\mbox{if  $|x|\leq c_a |y|$ and \,$l_{xy}\geq 2$, $l_{xy}\in\NN$, is such that $|y|\geq 32 a^{2l_{xy}}|x|$},\\[8pt]
C, \qquad \mbox{if}\quad c_a |y|\leq |x|\leq c_a^{-1}|y|.
\end{array}
\right.
\end{eqnarray}
\end{proposition}

\bp Let us focus first on the estimates for the second mixed
derivatives of $G$. The estimate for the case $c_a |y|\leq |x|\leq
c_a^{-1}|y|$ was proved in Theorem~\ref{t1.1}, and the bound for
$|x|\leq c_a |y|$ follows from the one for $|y|\leq c_a |x|$ by the
symmetry of Green's function. Hence, it is enough to consider the
case $|y|\leq c_a |x|$ only.

The function $x\mapsto \nabla_yG(x,y)$ is biharmonic in
$\Omega\setminus\{y\}$. We use Proposition~\ref{p7.2} with
$r=32|y|$ to write for $x\in \Omega\setminus B_{c_a^{-1}|y|}$
\begin{eqnarray}\label{eq7.21}
\nonumber |\nabla_x\nabla_y G(x,y)|^2 &  \leq &
\frac{C}{|x|^2\,|y|^3}\int_{C_{8|y|, 32|y|}}
|\nabla_y G(z,y)|^2\,dz\,\,\\[4pt]
&&\times\exp \left( -c\sum_{j=2}^{l}(|y|a^{2j}) {\rm
Cap}\,(\overline{C_{32|y|a^{2(j-1)},32|y|a^{2j}}}\setminus\Omega,\RR^3\setminus\{O\})\right),
\end{eqnarray}
where $l\geq 2$, $l\in\NN$, is such that $|x|\geq a^{2l}32|y|$.

Recall now
the function ${\cal R}$ introduced in the proof of
Proposition~\ref{p4.1}. If $y_0$ is a point on $\po$ such that
$|y-y_0|=d(y)$, then
\begin{equation}\label{eq7.22}
C_{8|y|,32|y|}\subset C_{6|y|,34|y|}(y_0),
\end{equation}

\noindent and  $\nabla_y G(z,y)={\cal R}(z,y)$ for every $z\in
C_{6|y|,34|y|}(y_0)$. Therefore,
\begin{eqnarray}\label{eq7.23}
&& \frac{1}{|x|^2\,|y|^3}\int_{C_{8|y|, 32|y|}} |\nabla_y
G(z,y)|^2\,dz\leq \frac{1}{|x|^2\,|y|^3}\int_{C_{6|y|,34|y|}(y_0)}
|{\cal R}(z,y)|^2\,dz\nonumber\\[4pt]
&&\qquad\leq \frac{C}{|x|^2\,d(y)^3}\int_{C_{3d(y)/2,6d(y)}(y_0)}
|{\cal R}(z,y)|^2\,dz\leq \frac{C}{|x|^2}\leq \frac{C}{|x-y|^2}.
\end{eqnarray}

\noindent The second inequality above follows from
Proposition~\ref{p3.4}, the third one has been proved in
(\ref{eq4.14})--(\ref{eq4.15}) and the last one owes to the
observation that
\begin{equation}\label{eq7.24}
|x-y|\leq |x|+|y|\leq (1+c_a)|x| \quad\mbox{whenever} \quad
|y|\leq c_a|x|.
\end{equation}

\noindent Combining (\ref{eq7.21})--(\ref{eq7.23}), we finish the
proof of the bound for the second mixed derivatives of Green's
function.

The proof of the estimate for $\nabla_y G$ follows a similar path,
and then the estimate for $\nabla_x G$ is a consequence of the
symmetry of Green's function. \ep

Analogously to (\ref{eq4.27})--(\ref{eq4.32}),
Proposition~\ref{p7.3} yields the following Corollary.

\begin{corollary}\label{c7.4} Suppose $u$ satisfies
\begin{equation}\label{eq7.25}
\Delta^2 u={\rm div}\,  f+h, \quad u\in \ring W_2^2(\Omega),
\end{equation}

\noindent for some functions $f=(f_1,f_2,f_3)$ and $h$. Fix some
$a\geq 4$ and let $c_a:=1/(32a^4)$. Then for any $x\in\Omega$
\begin{eqnarray}
\nonumber &&|\nabla u(x)|\\[4pt]\nonumber
&&\quad \leq C\int_{y\in\Omega:\,|y|\leq c_a |x|}
\exp \left( -c\sum_{j=2}^{l_{yx}}(|y|a^{2j}) {\rm
Cap}\,(\overline{C_{32|y|a^{2(j-1)},32|y|a^{2j}}}\setminus\Omega,\RR^3\setminus\{O\})\right)\,\nonumber\\[8pt]
&&\qquad\qquad\times\left(\frac{|f(y)|}{|x|}+|h(y)|\right)\,dy\nonumber\\[8pt]
&&\quad +\,C\int_{y\in\Omega:\,|x|\leq c_a |y|}
\exp \left( -c\sum_{j=2}^{l_{xy}}(|x|a^{2j}) {\rm
Cap}\,(\overline{C_{32|x|a^{2(j-1)},32|x|a^{2j}}}\setminus\Omega,\RR^3\setminus\{O\})\right)\,\nonumber\\[8pt]
&&\qquad\qquad\times\left(\frac{|f(y)|}{|y|}+|h(y)|\right)\,dy\nonumber \\[8pt]
&&\quad +\,C\int_{y\in\Omega:\,c_a |y|\leq |x|\leq c_a^{-1}|y|}
\left(\frac{|f(y)|}{|x-y|}+|h(y)|\right)\,dy,\nonumber
\end{eqnarray}
where in the first sum $l_{yx}\geq 2$, $l_{yx}\in\NN$, is such that $|x|\geq 32 a^{2l_{yx}}|y|$ and in the second sum $l_{xy}\geq 2$, $l_{xy}\in\NN$, is such that $|y|\geq 32 a^{2l_{xy}}|x|$.
\end{corollary}

\section{Necessary condition for 1-regularity}
\setcounter{equation}{0}

This section will be entirely devoted to the proof of the second
part of Theorem~\ref{t1.2}, i.e. the necessary condition for
1-regularity. We recall that ${\rm Cap}_P\,(K)={\rm
Cap}_P\,(K,\RR^3\setminus\{0\})$ for any compactum $K$ by
definition, and begin with

\vskip 0.08in

\noindent {\bf Step I: setup}. 
Suppose that for some $P\in\Pi_1$ the integral
in (\ref{eq1.10}) is convergent. For simplicity we shall assume that $a=2$. Any other value of $a$ could be treated in the exact same fashion.
Then for every $\eps>0$ there exists $N\in\NN$ such that
\begin{equation}\label{eq8.3}
\sum_{j=N}^\infty 2^{-j} \,{{\rm
Cap}_P\,(\overline{C_{2^{-j},2^{-j+2}}}\setminus\Omega,\RR^3\setminus
\{O\})}<\eps.
\end{equation}

Now let $K:=\overline{B_{2^{-N}}}\setminus \Omega$ and
$D:=\RR^3\setminus K$. We shall prove that the point $O$ is not
1-regular with respect to $D$, and therefore with respect to
$\Omega$, since $D$ coincides with $\Omega$ in a fixed
neighborhood of $O$ (see Proposition~\ref{p6.1}).

To this end, fix $P\in \Pi_1$ and let $\PP(x):=|x|P(x)$,
$x\in\RR^3$. Then take some cut-off $\eta\in C_0^\infty(B_{2})$
equal to 1 on $B_{3/2}$ and denote $f:=-[\Delta^2,\eta] \PP\in
C_0^\infty(B_2\setminus B_{3/2})$. Finally, let $V$ be a solution
of the boundary value problem
\begin{equation}\label{eq8.4}
\Delta^2 V=f \,\,{\mbox{in}}\,\,D, \quad V\in \ring W_2^2(D).
\end{equation}

\noindent Our goal is to show that $|\nabla V|$ does not vanish as
$x\to O$, $x\in D$.

Let us also consider the function $U:=V+\eta\,\PP$. One can check
that
\begin{equation}\label{eq8.5}
\Delta^2 U=0\,\,{\mbox{in}}\,\,D, \quad U=\PP
\,\,{\mbox{on}}\,\,K, \quad U\in \ring W_2^2(\RR^3).
\end{equation}

\noindent  Therefore, $U$ can be seen as a version of a biharmonic
potential. In fact, it is (\ref{eq8.5}) that gave an original idea
for the above definition of $V$.

\vskip 0.08in

\noindent {\bf Step II: main identity}. Let ${\cal B}$ denote the
bilinear form associated to the quadratic form in (\ref{eq2.3}),
i.e.
\begin{eqnarray}\label{eq8.6}
 {\cal B}(v,w)
&=&\int_{\RR}\int_{S^{2}}\Bigl[ (\delta_\omega v)(\delta_\omega
w)\,{\cal G}+2(\partial_t\nabla_\omega
v)\cdot(\partial_t\nabla_\omega
w)\,{\cal G}+ (\partial_t^2v)(\partial_t^2w)\,{\cal G} \nonumber\\[4pt]
&& -(\nabla_\omega v)\cdot (\nabla_\omega w) \Bigl(\partial_t^2
{\cal G}+\partial_t{\cal G}+2{\cal G}\Bigr)-(\partial_t
v)(\partial_t
w)\Bigl(2\partial_t^2{\cal G}+3\partial_t {\cal G}-{\cal G}\Bigr) \nonumber\\[4pt]
&&  +\frac 12\,v\,w\Bigl(\partial_t^4{\cal G}+2\partial_t^3{\cal
G} -\partial_t^2{\cal G}-2\partial_t{\cal G}\Bigr)\Bigr]\,d\omega
dt.
\end{eqnarray}

\noindent As before, we fix some point $\xi\in\RR^3$, $\tau:=\log
|\xi|^{-1}$ and let ${\cal G}(t)=g(t-\tau)$, $t\in\RR$. By ${\cal
B}_\tau (v,w)$ we shall denote ${\cal B}(v,w)$ for this particular
choice of ${\cal G}$. Then
\begin{eqnarray}\label{eq8.7}\nonumber
&& \int_{\RR^3}\Delta U(x)\,\Delta\bigg(\PP(x)|x|^{-1}\,g(\log
(|\xi|/|x|))\bigg)\,dx\\[4pt]
&&\qquad + \int_{\RR^3}\Delta
\PP(x)\,\Delta\bigg(U(x)|x|^{-1}\,g(\log (|\xi|/|x|))\bigg)\,dx=
2{\cal B}_{\tau} (u,q),
\end{eqnarray}

\noindent where $u=e^t(U\circ \varkappa^{-1})$ and $q=e^t(\PP\circ
\varkappa^{-1})=P\circ \varkappa^{-1}$.

The identity above can be proved along the lines of the argument
for Lemma~\ref{l2.1}, as long as the integration by parts and
absence of the boundary terms is justified. To this end, we note
that for any fixed $\xi\in\RR^3$ the function $x\mapsto g(\log
(|\xi|/|x|))$ is bounded by a constant as $|x|\to\infty$, while
$x\mapsto |x|^{-1}\,g(\log (|\xi|/|x|))$ is bounded by a constant
as $x\to O$. If $v_s\in C_0^\infty(D)$, $s\in\NN$, is a collection
of functions approximating $V$ in the $\ring W_2^2(D)$-norm, we
let $u_s:=v_s+\eta \PP$. Then $u_s$ converges to $U$ in $\ring
W_2^2(\RR^3)$. This, combined with the above observation about the
behavior of $g$, shows that it suffices to prove (\ref{eq8.7}) for
$u_s$ in place of $U$. Finally, since $u_s$ is compactly supported
in $\RR^3$ and is equal to $\PP$ in a neighborhood of 0, it is a
matter of direct calculation to establish (\ref{eq8.7}).

Since $(8\pi)^{-1}|x|$ is the fundamental solution of the
bilaplacian,
\begin{equation}\label{eq8.8}
\Delta^2
\PP(x)=\Delta^2(b_0|x|+b_1x_1+b_2x_2+b_3x_3)=(8\pi)^{-1}b_0\,\delta(x),
\end{equation}

\noindent where $\delta$ is the Dirac delta function. Therefore,
the second term on the left-hand side of (\ref{eq8.7}) is equal
(modulo a multiplicative constant) to $U(0)=0$.

Going further, we show that
\begin{equation}\label{eq8.9}
\int_{\RR^3}\Delta U(x)\,\Delta\bigg((U(x)-\PP(x))|x|^{-1}\,g(\log
(|\xi|/|x|))\bigg)\,dx=0.
\end{equation}

\noindent Indeed, the expression in (\ref{eq8.9}) is equal to
\begin{eqnarray}\label{eq8.10}\nonumber
&& \int_{\RR^3}\Delta U(x)\,\Delta\bigg(V(x)|x|^{-1}\,g(\log
(|\xi|/|x|))\bigg)\,dx\\[4pt]
&&\qquad + \int_{\RR^3}\Delta
U(x)\,\Delta\bigg((\eta(x)-1)\PP(x)|x|^{-1}\,g(\log
(|\xi|/|x|))\bigg)\,dx.
\end{eqnarray}

\noindent Then, using the aforementioned approximation by $v_s$,
$s\in\NN$, in the first integral, an observation that ${\rm
supp}\,(\eta-1) \PP\subset D$ in the second one, and the
biharmonicity of $U$ in $D$ we arrive at (\ref{eq8.9}).

Now the combination of (\ref{eq8.7})--(\ref{eq8.10}) leads to the
identity
\begin{equation}\label{eq8.11}
\int_{\RR^3}\Delta U(x)\,\Delta\bigg(U(x)|x|^{-1}\,g(\log
(|\xi|/|x|))\bigg)\,dx=2{\cal B}_\tau(u,q).
\end{equation}

\noindent Finally, since the identity (\ref{eq2.3}) holds for the
function $U$, (\ref{eq8.11}) implies that
\begin{equation}\label{eq8.12}
{\cal B}_{\tau}(u,u)=2{\cal B}_\tau(u,q).
\end{equation}

Recall now that $g$ is a fundamental solution of the equation
(\ref{eq2.7}), and therefore with the notation
\begin{eqnarray}\label{eq8.13}
&&\widetilde {\cal B}_\tau(v,w) =\int_{\RR}\int_{S^{2}}\Bigl[
(\delta_\omega v)(\delta_\omega
w)\,g(t-\tau)+2(\partial_t\nabla_\omega
v)\cdot(\partial_t\nabla_\omega
w)\,g(t-\tau) \nonumber\\[4pt]
&& \qquad+ (\partial_t^2v)(\partial_t^2w)\,g(t-\tau)
-(\nabla_\omega v)\cdot (\nabla_\omega w) \Bigl(\partial_t^2
g(t-\tau)+\partial_tg(t-\tau)+2g(t-\tau)\Bigr)\nonumber\\[4pt]
&&\qquad -(\partial_t v)(\partial_t
w)\Bigl(2\partial_t^2g(t-\tau)+3\partial_t
g(t-\tau)-g(t-\tau)\Bigr)\Bigr]\,d\omega dt,
\end{eqnarray}

\noindent we have
\begin{equation}\label{eq8.14}
{\cal B}_{\tau}(v,w)= \widetilde {\cal
B}_\tau(v,w)+\frac{1}{2}\int_{S^2}v(\tau,\omega)w(\tau,\omega)\,d\omega.
\end{equation}

\noindent Then the equality (\ref{eq8.12}) can be written as
\begin{equation}\label{eq8.15}
\int_{S^2}(u(\tau,\omega)-q(\tau,\omega))^2\,d\omega
=\int_{S^2}q^2(\tau,\omega)\,d\omega+4 \widetilde {\cal
B}_\tau(u,q)-2\widetilde {\cal B}_\tau(u,u),
\end{equation}

\noindent so that if $|\xi|<3/2$, $\tau=\log |\xi|^{-1}$,
\begin{equation}\label{eq8.16}
\int_{S^2}v^2(\tau,\omega)\,d\omega
=\int_{S^2}q^2(\tau,\omega)\,d\omega+4 \widetilde {\cal
B}_\tau(u,q)-2\widetilde {\cal B}_\tau(u,u),
\end{equation}

\noindent where $v=e^t(V\circ \varkappa^{-1})$.

\vskip 0.08in

 The identity (\ref{eq8.16}) is our major starting
point. We shall show that $\widetilde {\cal B}_\tau(u,q)$ and
$\widetilde {\cal B}_\tau(u,u)$ can be estimated in terms of the
series in (\ref{eq8.3}) and hence, can be made arbitrarily small
by shrinking $\eps$ in (\ref{eq8.3}). On the other hand,
\begin{equation}\label{eq8.17}
\int_{S^2}q^2(\tau,\omega)\,d\omega =
\int_{S^2}\Bigl(b_0^2+\sum_{i=1}^3 b_i^2
\omega_i^2\Bigr)\,d\omega=4\pi b_0^2+\frac{4\pi}{3}\sum_{i=1}^3
b_i^2,
\end{equation}

\noindent so that
\begin{equation}\label{eq8.18}
\frac{4\pi}{3}\leq \int_{S^2}q^2(\tau,\omega)\,d\omega \leq 4\pi.
\end{equation}

\noindent Therefore, by (\ref{eq8.16}),
\begin{equation}\label{eq8.19}
\int_{S^2}v^2(\tau,\omega)\,d\omega=\frac{C}{|\xi|^4}\int_{S_{|\xi|}}V^2(\xi)\,d\sigma_\xi
\end{equation}

\noindent  does not vanish as $\xi\to O$, which means that $\nabla
V$ does not vanish at $O$ either, as desired. It remains to
estimate $\widetilde {\cal B}_\tau(u,q)$ and $\widetilde {\cal
B}_\tau(u,u)$.

\vskip 0.08in

\noindent {\bf Step III: estimate for $\widetilde {\cal
B}_\tau(u,q)$}. Since $q=P\circ\varkappa^{-1}$ is independent of
$t$,
\begin{eqnarray}\label{eq8.20}
\widetilde {\cal B}_\tau(u,q) &=&\int_{\RR}\int_{S^{2}}\Bigl[
(\delta_\omega u)(\delta_\omega
q)\,g(t-\tau) \nonumber\\[4pt]
&& -(\nabla_\omega u)\cdot (\nabla_\omega q) \Bigl(\partial_t^2
g(t-\tau)+\partial_tg(t-\tau)+2g(t-\tau)\Bigr)\Bigr]\,d\omega dt.
\end{eqnarray}

\noindent Next, $\delta_\omega \omega_i=-2\omega_i$ for $i=1,2,3$,
and therefore $\delta_\omega q=-2\sum_{i=1}^3b_i\omega_i$, so that
\begin{eqnarray}\label{eq8.21}
\widetilde {\cal B}_\tau(u,q) &=&\int_{\RR}\int_{S^{2}}\Bigl[2b_0
\delta_\omega u\,g(t-\tau) -(\nabla_\omega u)\cdot (\nabla_\omega
q) \Bigl(\partial_t^2
g(t-\tau)+\partial_tg(t-\tau)\Bigr)\Bigr]\,d\omega dt\nonumber\\[4pt]
&=& - \int_{\RR}\int_{S^{2}}\Bigl[(\nabla_\omega u)\cdot
(\nabla_\omega q) \Bigl(\partial_t^2
g(t-\tau)+\partial_tg(t-\tau)\Bigr)\Bigr]\,d\omega dt
\nonumber\\[4pt]
&\leq &\left(\int_{\RR}\int_{S^{2}}\Bigl[|\nabla_\omega u|^2
\Bigl(-\partial_t^2
g(t-\tau)-\partial_tg(t-\tau)\Bigr)\Bigr]\,d\omega dt\right)^{1/2}\nonumber\\[4pt]
&& \times \left(\int_{S^{2}} |\nabla_\omega q|^2\,d\omega
\int_{\RR}\Bigl(-\partial_t^2
g(t-\tau)-\partial_tg(t-\tau)\Bigr)\,dt\right)^{1/2}=:I_1\times
I_2,
\end{eqnarray}

\noindent using the Cauchy-Schwarz inequality and the positivity
of the weight function (see (\ref{eq2.17})).

Inspecting the argument of Lemma~\ref{l2.3} one can see that
\begin{equation}\label{eq8.22}
I_1\leq (\widetilde {\cal B}_\tau(u,u))^{1/2}.
\end{equation}

\noindent On the other hand,
\begin{eqnarray}\label{eq8.23}\nonumber
I_2^2&=& \frac{8\pi}{3}\sum_{i=1}^3b_i^2
\int_{\RR}\Bigl(-\partial_t^2
g(t-\tau)-\partial_tg(t-\tau)\Bigr)\,dt\\[4pt]
&=& \frac{8\pi}{9}\sum_{i=1}^3b_i^2 \left(\int_{-\infty}^\tau
e^{t-\tau}\,dt+\int_\tau^\infty
e^{-2(t-\tau)}\,dt\right)=\frac{4\pi}{9}\sum_{i=1}^3b_i^2 \leq
\frac{4\pi}{9}.
\end{eqnarray}

\noindent Therefore,
\begin{equation}\label{eq8.24}
\widetilde {\cal B}_\tau(u,q)\leq \frac{2\sqrt \pi}{3}(\widetilde
{\cal B}_\tau(u,u))^{1/2}.
\end{equation}

\vskip 0.08in

\noindent {\bf Step IV: estimate for $\widetilde {\cal
B}_\tau(u,u)$, the setup}. Let us now focus on the estimate for
$\widetilde {\cal B}_\tau(u,u)$. To this end, consider the covering
of $K=\overline{B_{2^{-N}}}\setminus \Omega$ by the sets $K\cap
C_{2^{-j},2^{-j+2}}$, $j\geq N$, and observe that
\begin{eqnarray}\label{eq8.25}
&& K\cap {\overline{C_{2^{-j},2^{-j+2}}}}=
{\overline{C_{2^{-j},2^{-j+2}}}}\setminus \Omega, \quad j\geq
N+2,\\[4pt]\label{eq8.26}
&&  K\cap {\overline{C_{2^{-j},2^{-j+2}}}}\subseteq
{\overline{C_{2^{-j},2^{-j+2}}}}\setminus \Omega,\quad j=N, N+1.
\end{eqnarray}

\noindent Let $\{\eta^j\}_{j=N}^\infty$ be the corresponding
partition of unity such that
\begin{equation}\label{eq8.27}
\eta^j\in C_0^\infty(C_{2^{-j},2^{-j+2}}),\quad |\nabla^k
\eta^j|\leq C2^{kj}, \quad k=0,1,2,\quad\mbox{ and }\quad
\sum_{j=N}^\infty \eta^j=1.
\end{equation}

\noindent By $U^j$ we denote the capacitary potential of $K\cap
\overline{C_{2^{-j},2^{-j+2}}}$ with the boundary data $P$, i.e.
the minimizer for the optimization problem
\begin{eqnarray}\label{eq8.28}\nonumber
&&\inf\Bigg\{\int_{C_{2^{-j-2},2^{-j+4}}}(\Delta
u(x))^2\,dx:\,\,u\in \ring W_2^2(C_{2^{-j-2},2^{-j+4}}),\\[4pt]
&&\qquad\qquad\qquad\qquad\,u=P\mbox{ in a neighborhood of }K\cap
{\overline{C_{2^{-j},2^{-j+2}}}}\Bigg\}.
\end{eqnarray}

\noindent Such $U^j$ always exists and belongs to $\ring
W_2^2(C_{2^{-j-2},2^{-j+4}})$ since $P$ is an infinitely
differentiable function in a neighborhood of $K\cap
{\overline{C_{2^{-j},2^{-j+2}}}}$. The infimum above is equal to
\begin{equation}\label{eq8.29}
{\rm Cap}_P \{K\cap
{\overline{C_{2^{-j},2^{-j+2}}}},C_{2^{-j-2},2^{-j+4}}\}\approx
{\rm Cap}_P \{K\cap
{\overline{C_{2^{-j},2^{-j+2}}}},\RR^3\setminus\{O\}\}.
\end{equation}

Let us now define the function
\begin{equation}\label{eq8.30}
T(x):=\sum_{j=N}^\infty |x|\eta^j(x) U^j(x),\quad x\in \RR^3,
\end{equation}

\noindent and let $\vartheta:=e^t(T\circ \varkappa^{-1})$. Then by
the Cauchy-Schwarz inequality
\begin{equation}\label{eq8.31}
\widetilde {\cal B}_\tau(\vartheta,\vartheta)\leq C
\sum_{k=0}^2\sum_{j=N}^\infty
\int_{C_{2^{-j},2^{-j+2}}}\frac{|\nabla^k
(U^{j}(x))|^2}{|x|^{3-2k}}\,dx.
\end{equation}

\noindent Next, since $U^j\in \ring W_2^2(C_{2^{-j-2},2^{-j+4}})$,
the Hardy's inequality allows us to write
\begin{eqnarray}\label{eq8.32}\nonumber
&& \widetilde {\cal B}_\tau(\vartheta,\vartheta)\leq
C\sum_{j=N}^\infty 2^{-j}\int_{C_{2^{-j-2},2^{-j+4}}}|\nabla^2
U^{j}(x)|^2\,dx\leq C \sum_{j=N}^\infty
2^{-j}\int_{C_{2^{-j-2},2^{-j+4}}}|\Delta U^{j}(x)|^2\,dx\\[4pt]\nonumber
&&\quad\leq C\sum_{j=N}^\infty 2^{-j}\,{\rm Cap}_P \{K\cap
{\overline{C_{2^{-j},2^{-j+2}}}},\RR^3\setminus\{O\}\}\\[4pt]
&&\quad\leq C\sum_{j=N}^\infty 2^{-j}\,{\rm Cap}_P
\{{\overline{C_{2^{-j},2^{-j+2}}}}\setminus
\Omega,\RR^3\setminus\{O\}\} <C\eps,
\end{eqnarray}

\noindent by (\ref{eq8.29}), (\ref{eq8.25})--(\ref{eq8.26}), the
monotonicity property (\ref{eq5.4}), and (\ref{eq8.3}).

Having (\ref{eq8.32}) at hand, we need to consider the difference
$U-T$ in order to obtain the estimate for $\widetilde {\cal
B}_\tau(u,u)$. Let us denote $W:=U-T$, $w:=e^t(W\circ
\varkappa^{-1})$.

\vskip 0.08in

\noindent {\bf Step V: estimate for  ${\cal B}_\tau(w,w)$}. First
of all, one can show that $W\in \ring W_2^2(D)$. Indeed, both $U$
and $T$ belong to $\ring W_2^2(\RR^3)$. For $U$ this was pointed
out in (\ref{eq8.5}), the statement about $T$ can be proved along
the lines of (\ref{eq8.31})--(\ref{eq8.32}):
\begin{eqnarray}\label{eq8.32.1}\nonumber
&& \|T\|_{\ring W_2^2(\RR^3)}\leq C\sum_{k=0}^2\sum_{j=N}^\infty
2^{j(4-2k)}
\int_{C_{2^{-j},2^{-j+2}}}|\nabla^k (|x|U^{j}(x))|^2\,dx \\[4pt]\nonumber
&&\quad\leq C \sum_{j=N}^\infty
2^{-2j}\int_{C_{2^{-j-2},2^{-j+4}}}|\Delta U^{j}(x)|^2\,dx\\[4pt]
&&\quad\leq C\sum_{j=N}^\infty 2^{-2j}\,{\rm Cap}_P
\{{\overline{C_{2^{-j},2^{-j+2}}}}\setminus
\Omega,\RR^3\setminus\{O\}\} <C\eps.
\end{eqnarray}

\noindent In addition to (\ref{eq8.32.1}), we know that
$U={\mathbb P}$ on the boundary of $K$ by definition, and $\eta^j
U^j=U^j=P={\mathbb P}/|x|$ on the boundary of $K\cap
{\overline{C_{2^{-j},2^{-j+2}}}}$. Since by (\ref{eq8.30}) the
function $W$ is given by $\sum_{j=N}^\infty \eta^j(U-|x|U^j)$ in a
neighborhood of $K$, it vanishes on $\partial K$ (in the sense of
$\ring W_2^2(D)$).

Furthermore, $\Delta^2W=-\Delta^2 T$ in $D$ by (\ref{eq8.5}).
Then, with the notation $w:=e^t(W\circ \varkappa^{-1})$ we have
the formula
\begin{eqnarray}\label{eq8.33}\nonumber
 {\cal B}_\tau(w,w)&=&\int_{\RR^3}\Delta
W(x)\,\Delta\bigg(W(x)|x|^{-1}\,g(\log
(|\xi|/|x|))\bigg)\,dx\\[4pt]
& =& -\int_{\RR^3}\Delta T(x)\,\Delta\bigg(W(x)|x|^{-1}\,g(\log
(|\xi|/|x|))\bigg)\,dx.
\end{eqnarray}

\noindent
 In what follows we will show that
\begin{equation}\label{eq8.34}
-\int_{\RR^3}\Delta T(x)\,\Delta\bigg(W(x)|x|^{-1}\,g(\log
(|\xi|/|x|))\bigg)\,dx\leq C\eps^{1/2}({\cal B}_\tau(w,w))^{1/2}.
\end{equation}

Observe that according to (\ref{eq8.33}) and (\ref{eq2.11}) the
expression on the left-hand side of (\ref{eq8.34}) is positive.
Next, analogously to (\ref{eq2.5}),
\begin{eqnarray}\label{eq8.35}
&&\hskip -1.3cm -\int_{\RR^3}\Delta
T(x)\,\Delta\bigg(W(x)|x|^{-1}\,g(\log (|\xi|/|x|))\bigg)\,dx
\nonumber\\[4pt]
&&\hskip -1.3cm\quad = -\int_{\RR}\int_{S^2}
\Bigl(\partial_t^2\vartheta-3\partial_t\vartheta+2\vartheta+\delta_\omega
\vartheta\Bigr)\Bigl(g(t-\tau)\,\delta_\omega
w+g(t-\tau)\,\partial_t^2w\nonumber\\[4pt]
&&\quad\quad +(2\partial_t
g(t-\tau)-g(t-\tau))\,\partial_tw+(\partial_t^2g(t-\tau)-\partial_t
g(t-\tau))\,w\Bigr)\,d\omega dt.
\end{eqnarray}

\noindent Now recall the formula for $-(2\partial_t^2g+3\partial_t
g-g)$ from (\ref{eq2.18}). It is evident that for any
$D_t=\sum_{i=0}^4\alpha_i\partial_t^i$, $\alpha_i\in\RR$,  we have
\begin{equation}\label{eq8.36}
|D_tg|\leq C(-2\partial_t^2g-3\partial_t g+g),
\end{equation}

\noindent where $C$ generally depends on $D_t$, i.e. on
$\{\alpha_i\}_{i=0}^4$. Hence, for every such $D_t$
\begin{eqnarray}\label{eq8.37}
&&\hskip -1cm
\int_{\RR}\int_{S^2}(\partial_tw)^2\,|D_tg(t-\tau)|\,d\omega
dt\nonumber \\[4pt]
&&\hskip -1cm\quad \leq -\,C\int_{\RR}\int_{S^2}(\partial_tw)^2
(2\partial_t^2g(t-\tau)+3\partial_t g(t-\tau)-g(t-\tau))\,d\omega
dt\leq C\widetilde {\cal B}_\tau(w,w),
\end{eqnarray}

\noindent where the last inequality follows from the calculations
in Lemma~\ref{l2.3}. Then, using (\ref{eq8.31})--(\ref{eq8.32}),
we have
\begin{eqnarray}\label{eq8.38}
&&\int_{\RR}\int_{S^2}|\partial_t^k\nabla_\omega^i\vartheta|\,
|\partial_tw|\,|D_tg(t-\tau)|\,d\omega
dt\nonumber \\[4pt]
&& \qquad \leq C
\left(\int_{\RR}\int_{S^2}|\partial_t^k\nabla_\omega^i\vartheta|^2\,
|D_tg(t-\tau)|\,d\omega dt \right)^{1/2} (\widetilde {\cal
B}_\tau(w,w))^{1/2}\nonumber \\[4pt]
&& \qquad \leq C \sum_{j=0}^2
\left(\int_{\RR^3}\frac{|\nabla^jT(x)|^2}{|x|^{5-2j}}\,dx
\right)^{1/2} (\widetilde {\cal B}_\tau(w,w))^{1/2}\leq
C\eps^{1/2} (\widetilde {\cal B}_\tau(w,w))^{1/2},
\end{eqnarray}

\noindent for $0\leq i+k\leq 2$.

For similar reasons,
\begin{equation}\label{eq8.39}
\int_{\RR}\int_{S^2}|\partial_t^k\nabla_\omega^i\vartheta|\,
|\partial_t^2w|\,g(t-\tau)\,d\omega dt\leq C\eps^{1/2} (\widetilde
{\cal B}_\tau(w,w))^{1/2},
\end{equation}

\noindent and
\begin{equation}\label{eq8.40}
\int_{\RR}\int_{S^2}|\partial_t^k\nabla_\omega^i\vartheta|\,
|\partial_t\nabla_\omega w|\,g(t-\tau)\,d\omega dt\leq C\eps^{1/2}
(\widetilde {\cal B}_\tau(w,w))^{1/2},
\end{equation}

\noindent for $0\leq i+k\leq 2$.

Invoking (\ref{eq8.38})--(\ref{eq8.40}) and integrating by parts,
we see that the expression in (\ref{eq8.35}) is bounded by
\begin{eqnarray}\label{eq8.41}
&&\Bigg|\int_{\RR}\int_{S^2}\Big(\delta_\omega\vartheta\delta_\omega
w\, g(t-\tau)-\nabla_\omega\vartheta\cdot\nabla_\omega
w\,(2\partial_t^2g(t-\tau)+2\partial_tg(t-\tau)+2g(t-\tau))\nonumber \\[4pt]
&&\qquad+\vartheta w
(\partial_t^4g(t-\tau)+2\partial_t^3g(t-\tau)-\partial_t^2g(t-\tau)-2\partial_tg(t-\tau))\Big)\,dtd\omega\Bigg|
\nonumber \\[4pt]
&&\qquad +C\eps^{1/2} (\widetilde {\cal B}_\tau(w,w))^{1/2}.
\end{eqnarray}

\noindent Also,
\begin{eqnarray}\label{eq8.42}
&&  \left|\int_\RR \int_{S^2} \left( \delta_\omega \vartheta \cdot
\delta_\omega w - 2 \nabla_\omega \vartheta \cdot \nabla_\omega w
\right)\, g \, dt\, d\omega \right|\nonumber\\[4pt]
&&\qquad \leq \left( \int_\RR \int_{S^2} \left[ (\delta_\omega
\vartheta)^2-2(\nabla_\omega \vartheta)^2\right]\,g\, dt\,d\omega
\right)^{1/2} \left( \int_\RR \int_{S^2} \left[ (\delta_\omega
w)^2-2(\nabla_\omega w)^2\right]\,g\, dt\,d\omega \right)^{1/2}
\nonumber\\[4pt]
&&\qquad\leq C \eps^{1/2} (\widetilde {\cal B}_\tau(w,w))^{1/2},
\end{eqnarray}

\noindent using (\ref{eq2.12}) and the Cauchy-Schwarz inequality
for the bilinear form on the left-hand side of (\ref{eq8.42}). In
view of (\ref{eq8.42}) and (\ref{eq2.7}) the expression in
(\ref{eq8.41}) is controlled by
\begin{eqnarray}\label{eq8.43}
&&\Bigg|\int_{\RR}\int_{S^2}\nabla_\omega\vartheta\cdot\nabla_\omega
w\,(-2\partial_t^2g(t-\tau)-2\partial_tg(t-\tau))\,dtd\omega\Bigg|
\nonumber \\[4pt]
&&\qquad +\frac 12\left|\int_{S^2}\vartheta(\tau,\omega)
w(\tau,\omega)\,d\omega\right|+C\eps^{1/2} (\widetilde {\cal
B}_\tau(w,w))^{1/2}\nonumber \\[4pt]
&&\leq C\eps^{1/2} (\widetilde {\cal B}_\tau(w,w))^{1/2}+\frac
12\left(\int_{S^2}\vartheta^2(\tau,\omega)\,d\omega\right)^{1/2}({\cal
B}_\tau(w,w))^{1/2}.
\end{eqnarray}

\noindent Here we used the positivity of
$-2\partial_t^2g-2\partial_tg$ (see (\ref{eq2.17})) and the
argument similar to (\ref{eq8.37})--(\ref{eq8.38}) to estimate the
first term. The bound for the second one follows from the
Cauchy-Schwarz inequality and (\ref{eq2.11}).

Finally, we claim that
\begin{equation}\label{eq8.44}
\int_{S^2}\vartheta^2(\tau,\omega)\,d\omega<C\eps.
\end{equation}

\noindent Indeed, by definition (\ref{eq8.44}) is equal to
\begin{eqnarray}\label{eq8.45}
&&\frac{1}{|\xi|^4}\int_{S_{|\xi|}}T^2(\xi)\,d\sigma_\xi\leq C
\sum_{j:\,2^{-j}\leq |\xi|\leq 2^{-j+2}}
\frac{1}{|\xi|^2}\int_{S_{|\xi|}}(U^j(\xi))^2\,d\sigma_\xi\nonumber\\[4pt]
&& \qquad  \leq C \sum_{j:\,2^{-j}\leq |\xi|\leq 2^{-j+2}}
\int_{\RR^3}\Delta \bigg(|x|U^j(x)\bigg)\Delta\bigg(U^j(x)\,g(\log
(|\xi|/|x|))\bigg)\,dx,
\end{eqnarray}

\noindent using (\ref{eq2.11}) for the function $x\mapsto
|x|\,U^j(x)$ in $\ring W_2^2(C_{2^{-j-2},2^{-j+4}})$. Finally, the
right-hand side of (\ref{eq8.45}) is bounded by
\begin{equation}\label{eq8.46}
C \sum_{j:\,2^{-j}\leq |\xi|\leq 2^{-j+2}}\sum_{k=0}^2
\int_{C_{2^{-j},2^{-j+2}}}\frac{|\nabla^k
(|x|U^{j}(x))|^2}{|x|^{5-2k}}\,dx<C\eps,
\end{equation}

\noindent by the estimate following (\ref{eq8.31}). This completes
the proof of (\ref{eq8.34}), which together with (\ref{eq8.33})
yields ${\cal B}_\tau(w,w)<\eps$.  and therefore,
\begin{equation}\label{eq8.47}
\widetilde {\cal B}_\tau(w,w)<  {\cal B}_\tau(w,w)<C\eps.
\end{equation}

\noindent The last estimate, in turn, implies that ${\cal
B}_\tau(u,u)<C\eps$ owing to the results of Step IV. At last, the
combination with (\ref{eq8.24}) finishes the argument. \ep

\section{Examples and further properties of ${\rm Cap}_P$ and ${\rm Cap}$.}
\setcounter{equation}{0}

\begin{lemma}\label{l9.1} Consider a domain $\Omega$ shaped as an
exterior of a cusp in some neighborhood of $O\in\po$, i.e.
\begin{equation}\label{eq9.1}
\Omega\cap
B_c=\{(r,\theta,\phi):\,0<r<c,\,\theta>h(r)\},\quad\mbox{for
some}\quad c>0,
\end{equation}

\noindent where $(r,\theta,\phi)$, $r\in(0,c)$,
$\theta\in[0,\pi]$, $\phi\in [0,2\pi)$, are  spherical coordinates
in $\RR^3$ and $h(r):(0,c)\to\RR$ is a nondecreasing function
satisfying the condition $h(br)\leq Ch(r)$ for some $b>1$ and all
$r\in (0,c)$.

Then
\begin{equation}\label{eq9.2}
O \,\,\mbox {is 1-regular} \quad\mbox{ if and only if }\quad
\int_{0}^{c} s^{-1}h(s)^2\,ds=+\infty.
\end{equation}
\end{lemma}

\bp We claim that for every $P\in\Pi_1$ and every $a\geq 4$
\begin{equation}\label{eq9.3}
{{\rm Cap}_P\,(\overline{C_{s,as}}\setminus\Omega,\RR^3\setminus
\{O\})}\geq  C s^{-1}h(s)^2, \qquad 0<s<c/a.
\end{equation}

Indeed, recall from Lemma~\ref{l5.2} that
\begin{equation}\label{eq9.4}
{{\rm Cap}_P\,(\overline{C_{s,as}}\setminus\Omega,\RR^3\setminus
\{O\})}\approx {{\rm
Cap}_P\,(\overline{C_{s,as}}\setminus\Omega,C_{s/2,2as})}.
\end{equation}

\noindent By definition of the capacity ${\rm Cap}_P$, for every
$\eps>0$ there exists some $u\in \ring W_2^2(C_{s/2,2as})$ such
that
\begin{equation}\label{eq9.5}
{{\rm
Cap}_P\,(\overline{C_{s,as}}\setminus\Omega,C_{s/2,2as})}+\eps\geq
C \int_{C_{s/2,2as}}(\Delta u(x))^2\,dx,
\end{equation}

\noindent and $u=P$ in a neighborhood of
$\overline{C_{s,as}}\setminus\Omega$. Since $u\in \ring
W_2^2(C_{s/2,2as})$, by Hardy's inequality
\begin{eqnarray}\label{eq9.6}
&& \int_{C_{s/2,2as}}(\Delta
u(x))^2\,dx=\int_{C_{s/2,2as}}|\nabla^2
u(x)|^2\,dx\nonumber\\[4pt]
&& \qquad \geq C \int_{C_{s/2,2as}}\Bigg(\frac {|u(x)|^2}{|x|^4}+\frac {|\nabla u(x)|^2}{|x|^2}\Bigg)\,dx
\geq C \int_{C_{s,as}\setminus\Omega}\Bigg(\frac {|P(x)|^2}{|x|^4}+\frac {|\nabla P(x)|^2}{|x|^2}\Bigg)\,dx\nonumber\\[4pt]
&& \qquad \geq C \int_{C_{s,as}\setminus\Omega}\Bigg(\frac {|P(x)|^2}{|x|^4}+\frac {|\nabla \left(|x|\,P(x)\right)|^2}{|x|^4}\Bigg)\,dx.
\end{eqnarray}

\noindent The contribution from $\frac {|P(x)|^2}{|x|^4}$ amounts to
\begin{eqnarray}\label{eq9.7}\nonumber
&&\hskip -1cm \frac{C}{s^4}\int_s^{as}\int_{0}^{h(r)}\int_0^{2\pi}\left(b_0+b_1\sin
\theta\cos\phi+b_2\sin\theta\sin\phi+b_3\cos\theta\right)^2\,sin\theta\,
r^2\,d\phi d\theta dr\\[4pt]
&&\hskip -1cm \quad \geq \frac{C}{s}\int_{0}^{h(s)}\left(b_0^2+
b_3^2\cos^2\theta+2b_0b_3\cos\theta\right)\,sin\theta\, d\theta\nonumber \\[4pt]
&&\hskip -1cm \quad \geq \frac{C}{s}\left(\cos\theta\Bigl(b_0^2
+\frac{b_3^2}{3}\cos^2\theta+b_0b_3\cos\theta\Bigr)
\right)\Bigg|_{h(s)}^0
\nonumber \\[4pt]
&&\hskip -1cm \quad \geq \frac{C}{s}\left(\cos\theta\Bigl(\frac 14\, b_0^2
+\Bigl(\frac {\sqrt 3} 2 \,b_0-\frac{\cos\theta}{\sqrt 3}\,b_3\Bigr)^2\Bigr)
\right)\Bigg|_{h(s)}^0 \geq \frac{C}{s}\,b_0^2
\cos\theta\,\Bigg|_{h(s)}^0\geq \frac{C}{s}\,b_0^2h(s)^2.
\end{eqnarray}

\noindent On the other hand,
\begin{equation}\label{eq9.7.1}
|\nabla (|x|\,P(x))|^2=\sum_{i=1}^3\left(b_0\frac{x_i}{|x|}+b_i\right)^2
\end{equation}

\noindent and for every $i=1,2,3$
\begin{equation}\label{eq9.7.2}
\left(b_0\frac{x_i}{|x|}+b_i\right)^2+b_0^2\approx b_i^2+b_0^2.
\end{equation}

\noindent Hence,

\begin{equation}\label{eq9.7.3}
\int_{C_{s,as}\setminus\Omega}\Bigg(\frac {|P(x)|^2}{|x|^4}+\frac {|\nabla \left(|x|\,P(x)\right)|^2}{|x|^4}\Bigg)\,dx\geq \frac{C}{s}\, h(s)^2 \sum_{i=0}^3b_i^2\geq \frac{C}{s}\, h(s)^2.
\end{equation}

\noindent
Now one can combine (\ref{eq9.5}), (\ref{eq9.6}), (\ref{eq9.7.3}) and let
$\eps\to 0$ to obtain (\ref{eq9.3}).

Conversely, we claim that there exists $P\in\Pi_1$ such that for
every $s\in(0,c/a)$
\begin{equation}\label{eq9.9}
{{\rm Cap}_P\,(\overline{C_{s,as}}\setminus\Omega,\RR^3\setminus
\{O\})}\leq  C s^{-1}h(s)^2.
\end{equation}

\noindent Indeed, let us take
\begin{equation}\label{eq9.10}
P(x):=\frac 12\Bigl(1-\frac{x_3}{|x|}\Bigr), \qquad x\in\RR^3.
\end{equation}

\noindent  Clearly, $P\in\Pi_1$. Next, we choose a function $U\in
\ring W_2^2(C_{s/2,2as})$ that is given by $P$ in a neighborhood
of $C_{s,as}\setminus\Omega$. To do this, let us introduce two
cut-off functions, $\zeta^\theta$ and $\zeta^r$, such that
\begin{equation}\label{eq9.11}
\zeta^{\theta}\in C_0^\infty\left(-1/2,2\right), \,\,
\zeta^\theta=1\,\,\mbox{on}\,\,[0,3/2]; \quad \zeta^r \in
C_0^\infty\left(1/2,2a\right), \,\,
\zeta^r=1\,\,\mbox{on}\,\,[3/4,3a/2].
\end{equation}

\noindent Then let
\begin{equation}\label{eq9.12}
u(r,\phi,\theta):=\frac
12(1-\cos\theta)\,\zeta^\theta\Bigl(\frac\theta{h(as)}\Bigr)\,\zeta^r\Bigl(\frac
rs\Bigr),
\end{equation}

\noindent so that
\begin{equation}\label{eq9.13}
u(r,\phi,\theta)=1\quad\mbox{whenever}\quad 0\leq \theta\leq
\frac{3h(as)}{2}\quad\mbox{and}\quad \frac{3s}{4}\leq r\leq
\frac{3as}{2},
\end{equation}

\noindent and
\begin{equation}\label{eq9.14}
u(r,\phi,\theta)=0\quad\mbox{whenever}\quad 2h(as)\leq \theta\leq
\pi\quad\mbox{or}\quad r\not\in \Bigl(\frac s2, 2as\Bigr).
\end{equation}

\noindent Finally, let $U:=u\circ \kappa$, where $\kappa$ is the
change of coordinates in (\ref{eq2.2}). Then
\begin{eqnarray}\nonumber
\int_{C_{s/2,2as}}|\Delta U(x)|^2\,dx =
C\int_{s/2}^{2as}\int_0^{2h(as)} \Bigg| \frac {1}{r^2}\,\partial_r
(r^2\partial_r
u)+\frac{1}{r^2\sin\theta}\,\partial_\theta(\sin\theta\,\partial_\theta
u)\Bigg|^2\sin\theta \,d\theta\,r^2\,dr,
\end{eqnarray}

\noindent since $u$ is independent of $\phi$. A straightforward
calculation shows that for $r$ and $\theta$ as above
\begin{equation}\label{eq9.15}
\Bigg| \frac {1}{r^2}\,\partial_r (r^2\partial_r
u)+\frac{1}{r^2\sin\theta}\,\partial_\theta(\sin\theta\,\partial_\theta
u)\Bigg|\leq \frac{C}{s^2},
\end{equation}

\noindent and therefore,
\begin{equation}\label{eq9.16}
\int_{C_{s/2,2as}}|\Delta U(x)|^2\,dx\leq C s^{-1}h(as)^2 \leq C
s^{-1}h(s)^2.
\end{equation}

\noindent If $a\leq b$, the last inequality follows from the fact
that $h$ is nondecreasing. If $a>b$, we have $h(ar)\leq C^m
h\Bigl(ab^{-m}r\Bigr)\leq C^{m+1}h(r)$ for
 $m\geq \log_b a-1$.
 
Combining (\ref{eq9.3}) and (\ref{eq9.9}) and employing 
Theorem~\ref{t1.2} together with the integral test for series convergence, we arrive at 
(\ref{eq9.2}). \ep

In order to state the next result, let us recall one of the
definitions of the harmonic capacity of a compact set. For an open
set $\Omega\subset\RR^3\setminus \{O\}$ and a compactum $e\subset
\Omega$
\begin{equation}\label{eq9.17}
{\rm cap}\,(e,\Omega):=\inf\Bigg\{\int_\Omega (\nabla
u(x))^2\,dx:\,\,u\in \ring W^1_2(\Omega),\,u=1\mbox{ in a
neighborhood of }e\Bigg\},
\end{equation}

\noindent is a harmonic capacity of the set $e$ relative to
$\Omega$. If $\Omega=\RR^3\setminus\{0\}$ then (\ref{eq9.17})
coincides with (\ref{eq1.5}).

\begin{lemma}\label{l9.2} Let $K$ be a compactum situated on the
set
\begin{equation}\label{eq9.18}
\{x\in\RR^3:\,b_0 |x|+b_1x_1+b_2x_2+b_3x_3=0\},\quad
b_i\in\RR,\,i=0,1,2,3,
\end{equation}

\noindent such that $O\not\in K$. If the harmonic capacity of $K$
equals zero, then
\begin{equation}\label{eq9.19}
{\rm Cap}_P(K,\RR^3\setminus\{0\})=0
\end{equation}

\noindent for
\begin{equation}\label{eq9.20}
P(x)=\frac{1}{\sqrt{b_0^2+b_1^2+b_2^2+b_3^2}}\,\Bigl(b_0+b_1\,\frac{x_1}{|x|}+b_2\,\frac{x_2}{|x|}
+b_3\,\frac{x_3}{|x|}\Bigr), \qquad x\in\RR^3\setminus \{O\}.
\end{equation}

\noindent In particular, ${\rm Cap}(K,\RR^3\setminus\{0\})=0$.
\end{lemma}

\bp By current assumptions, $O\not\in K$. Therefore, there exist
$s>0$ and $a>1$ such that $K\subset\overline{C_{s,as}}$. In the
course of proof some constants will depend on $s$ and $a$. That,
however, does not influence the result.

Since
\begin{equation}\label{eq9.21}
{\rm cap}(K, C_{s/2,2as})\approx {\rm
cap}(K,\RR^3\setminus\{0\})=0,
\end{equation}

\noindent for every $\eps>0$ there exists a compactum $K_\eps$
with a smooth boundary contained in the set (\ref{eq9.18}) and
such that
\begin{equation}\label{eq9.22}
K\subset K_\eps\subset C_{s/2,2as}\quad\mbox{and}\quad {\rm
cap}(K_\eps, C_{s/2,2as})<\eps.
\end{equation}

\noindent Let $u$ denote the harmonic potential of $K_\eps$, so
that
\begin{equation}\label{eq9.23}
u\in \ring W_2^1(C_{s/2,2as}),\quad u=1\mbox{ in
}K_\eps,\quad\Delta u=0 \mbox{ in }\RR^3\setminus K_\eps,\quad
\int_{C_{s/2,2as}}|\nabla u(x)|^2\,dx<\eps.
\end{equation}

\noindent Next, given $\alpha<1$ let
\begin{equation}\label{eq9.24}
v_\alpha(x)=\left\{\begin{array}{l}
\alpha^{-4}P(x)u^2(x)(2\alpha - u(x))^2,\qquad \mbox{if}\quad u(x)\leq\alpha,\\[4pt]
P(x),\qquad\qquad\qquad\qquad\qquad\quad\,\, \mbox{if}\quad
u(x)>\alpha,
\end{array}
\right.
\end{equation}

\noindent where $x\in C_{s/2,2as}$ and $P$ is defined by
(\ref{eq9.20}). Then $v_\alpha\in \ring W_2^2(C_{s/2,2as})$ by
(\ref{eq9.23}) and $v_\alpha=P$ in a neighborhood of $K$.
Therefore,
\begin{eqnarray}\label{eq9.25}\nonumber
&&{\rm Cap}_P(K,\RR^3\setminus\{0\})\approx {\rm
Cap}_P(K,C_{s/2,2as})\leq \int_{C_{s/2,2as}}|\Delta
v_\alpha(x)|^2\,dx\\[4pt]
&&\quad = \alpha^{-8}\int_{x:\,u(x)\leq\alpha}\Bigl|\Delta
\Bigl(P(x)u^2(x)(2\alpha -
u(x))^2\Bigr)\Bigr|^2\,dx+\int_{x:\,u(x)>\alpha}|\Delta
P(x)|^2\,dx.
\end{eqnarray}

\noindent We take $\alpha=\alpha(\eps)<1$ (close to 1) such that
the last term above is less than $\eps$. In addition, on the set
$\{x:\,u(x)\leq\alpha\}$
\begin{eqnarray}\label{eq9.26}\nonumber
&&\Bigl|\Delta\Bigl(u^2(x)(2\alpha - u(x))^2\Bigr)\Bigr|\leq
C|\nabla u|^2,\qquad\Bigl|\nabla P\cdot\nabla\Bigl(u^2(x)(2\alpha
-
u(x))^2\Bigr)\Bigr|\leq C|\nabla u|,\\[4pt]
&&\qquad\qquad\qquad \Bigl|\Delta P\Bigl(u^2(x)(2\alpha -
u(x))^2\Bigr)\Bigr|\leq C|u|,
\end{eqnarray}

\noindent so that
\begin{equation}\label{eq9.27}
\int_{x:\,u(x)\leq\alpha}\Bigl|\Delta \Bigl(P(x)u^2(x)(2\alpha -
u(x))^2\Bigr)\Bigr|^2\,dx\leq
C\eps+C\int_{x:\,u(x)\leq\alpha}|P(x)|^2 |\nabla u|^4\,dx,
\end{equation}

\noindent by (\ref{eq9.23}).

It remains to estimate the last integral above. Let us denote by
$\{B_i\}_{i=1}^\infty$ a Whitney decomposition of the set
$C_{s/2,2as}\setminus K_\eps$, i.e. a collection of balls such
that
\begin{equation}\label{eq9.28}
\bigcup_{i=1}^\infty B_i=C_{s/2,2as}\setminus K_\eps,\quad
\sum_{i=1}^\infty \chi_{B_i}\leq C,\quad r(B_i)\approx {\rm
dist}\,\Bigl(B_i,\partial(C_{s/2,2as}\setminus K_\eps)\Bigr),
\end{equation}

\noindent where $r(B_i)$ denotes the radius of $B_i$. Observe that
\begin{eqnarray}\nonumber
&& |u(x)|\leq 1,\quad |P(x)|\leq Cr_i,\quad\mbox{if}\,\, x\in B_i
\,\,\mbox{such that}\,\, {\rm dist}\,(B_i,\partial
C_{s/2,2as})\geq {\rm dist}\,(B_i, K_\eps),\\[4pt]\nonumber
&& |u(x)|\leq Cr_i,\quad |P(x)|\leq C,\quad\mbox{if}\,\, x\in B_i
\,\,\mbox{such that}\,\, {\rm dist}\,(B_i,\partial
C_{s/2,2as})\leq {\rm dist}\,(B_i, K_\eps).
\end{eqnarray}

\noindent Since $u$ is harmonic in $C_{s/2,2as}\setminus K_\eps$,
\begin{equation}\label{eq9.29}
|\nabla u|^2\leq \frac{C}{r_i^5}\int_{B_{i}}|u(x)|^2\,dx.
\end{equation}

\noindent Therefore, $|P||\nabla u|\leq C$ on
$C_{s/2,2as}\setminus K_\eps$ and
\begin{equation}\label{eq9.30}
\int_{C_{s/2,2as}}|P(x)|^2 |\nabla u|^4\,dx\leq
\int_{C_{s/2,2as}}|\nabla u|^2\,dx<\eps.
\end{equation}

\noindent Letting $\eps\to 0$, we finish the argument.\ep
\begin{corollary}\label{c9.3}
Let $\Omega$ be a domain in $\RR^3$ such that $O\in\po$ and the
complement of $\Omega$ is a compactum of zero harmonic capacity
situated on the set (\ref{eq9.18}). Then the point $O$ is not
1-regular.
\end{corollary}

\bp By Lemma~\ref{l9.2} for the choice of $P$ in (\ref{eq9.20})
\begin{equation}\label{eq9.31}
{{\rm Cap}_P\,(\overline{C_{s,as}}\setminus\Omega,\RR^3\setminus
\{O\})}\,ds=0,
\end{equation}

\noindent for every $s>0$, $a>1$. One can see that such $P$ does
not depend on $s$ and $a$, but only on the initial cone containing
the complement of $\Omega$. Therefore,
\begin{equation}\label{eq9.32}
\inf_{P\in\Pi_1}\sum_{j=0}^\infty a^{-j} \,{{\rm
Cap}_P\,(\overline{C_{a^{-j},a^{-j+2}}}\setminus\Omega)} =0,
\end{equation}

\noindent and hence $O$ is not 1-regular by Theorem~\ref{t1.2}.\ep

\vskip 0.08in

\noindent{\it Remark}.\,The set defined by (\ref{eq9.18}) is
either a circular cone with the vertex at $O$ or a plane
containing $O$. Indeed, the set (\ref{eq9.18}) is formed by the
rays originating at $O$ and passing through the intersection of
the plane $b_0+b_1x_1+b_2x_2+b_3x_3=0$ with the unit sphere. If
this plane passes through the origin ($b_0=0$), it is actually the
set (\ref{eq9.18}). If it does not, then its intersection with
$S^2$ is a circle giving rise to the corresponding circular cone.

Due to the particular form of elements in the space $\Pi_1$ such
sets play a special role for our concept of the capacity and for
1-regularity. This observation is, in particular, supported by
Lemma~\ref{l9.2} and the upcoming example.

\vskip 0.08in

We consider a domain whose complement consists of a set of points
such that in each dyadic spherical layer three of the points
belong to a fixed circular cone, while the fourth one does not.
The result below shows that in this case the origin is 1-regular
provided the deviation of the fourth point is large enough in a
certain sense. The details are as follows.

\begin{lemma}\label{l9.4} Fix some $a\geq 4$.
Consider a domain $\Omega$ such that in
some neighborhood of the origin its complement consists of the set
of points
\begin{equation}\label{eq9.33}
\bigcup_{k}
\left\{A_1^k=(a^{-k},0,\alpha),\,\,A_2^k=(a^{-k},\pi/2,\alpha),\,\,
A_3^k=(a^{-k},\pi,\alpha),\,\,A_4^k=(a^{-k+1/2},3\pi/2,\beta_k)\right\},
\end{equation}

\noindent where the points are represented in spherical
coordinates $(r,\phi,\theta)$, $r\in(0,c)$ for some $c>0$,
$\theta\in[0,\pi]$, $\phi\in [0,2\pi)$, $k\in\NN\cap(1/2-\log_a c,\infty)$. Assume, in addition, that
\begin{equation}\label{eq9.34}
0<\alpha<\pi/2,\quad 0\leq |\beta_k-\alpha|<\alpha/2, \quad
\forall\,\,k\in\NN\cap(1/2-\log_a c,\infty).
\end{equation}

\noindent Then
\begin{equation}\label{eq9.35}
\sum_{k} a^{-k} \,\inf_{P\in\Pi_1}{{\rm
Cap}_P\,(\overline{C_{a^{-k},a^{-k+1}}}\setminus\Omega)} \geq C\sum_{k} (\beta_k-\alpha)^2,
\end{equation}

\noindent where $C=C(\alpha)>0$ and the summation is over $k\in\NN\cap(1/2-\log_a c,\infty)$. In particular,
\begin{equation}\label{eq9.36}
{\mbox if}\quad \sum_{k}
(\beta_k-\alpha)^2=+\infty\quad\mbox{then}\,\, O \,\,\mbox {is
1-regular}.
\end{equation}
\end{lemma}

\bp To begin, let us observe that
 in the spherical
layer $\overline{C_{a^{-k},a^{-k+1}}}$ there are exactly four points that
belong to the complement of $\Omega$, namely, $A_i^k$,
$i=1,2,3,4$. We aim to show that 
\begin{equation}\label{eq9.38}
{{\rm Cap}\,(\overline{C_{a^{-k},a^{-k+1}}}\setminus\Omega,\RR^3\setminus
\{O\})}\geq Ca^k (\beta_k-\alpha)^2.
\end{equation}

Fix $s=a^{-k}$. Take some $P\in\Pi_1$ and consider the distribution
\begin{equation}\label{eq9.39}
T^k(x):=\sum_{i=1}^4 P(A_i^k)\delta (x-A_i^k).
\end{equation}

\noindent Then for every $u\in \ring W_2^2(C_{s/2,2as})$ such that
$u=P$ in a neighborhood of $\{A_i^k, i=1,2,3,4\}$, we have
\begin{equation}\label{eq9.40}
\langle T^k,P\rangle=\sum_{i=1}^4 P(A_i^k)^2.
\end{equation}

On the other hand, since $T^k$ is supported in the set $\{A_i^k,
i=1,2,3,4\}$,
\begin{equation}\label{eq9.41}
\langle T^k,P\rangle= -\langle \Delta E \ast T^k,u\rangle=
-\langle E \ast T^k, \Delta u\rangle,
\end{equation}

\noindent where $E(x)={1}/({4\pi|x|})$ is the fundamental solution
for the Laplacian. By the Cauchy-Schwarz inequality
\begin{eqnarray}\nonumber
|\langle T^k,P\rangle|^2&\leq &\|E\ast T^k\|_{L^2(C_{s/2,2as})}^2
\|\Delta u\|_{L^2(C_{s/2,2as})}^2\\[4pt]\label{eq9.42}
&\leq & C s\,\sum_{i=1}^4 P(A_i^k)^2\, {{\rm
Cap}_P\,(\overline{C_{s,as}}\setminus\Omega, C_{s/2,2as})}.
\end{eqnarray}

\noindent Therefore, combining (\ref{eq9.40})--(\ref{eq9.42}) and
taking the infimum in $P$, we obtain the estimate
\begin{equation}\label{eq9.43}
{{\rm Cap}\,(\overline{C_{s,as}}\setminus\Omega,\RR^3\setminus
\{O\})}\geq Ca^k \inf_{P\in\Pi_1}\sum_{i=1}^4
P(A_i^k)^2=Ca^k\inf_{b\in\RR^4:\,\|b\|=1} b\,M M^\bot\, b^\bot,
\end{equation}

\noindent where $b=(b_0, b_1, b_2, b_3)$,
\begin{equation}\label{eq9.44}
M=\left(
\begin{array}{cccc}
1 & 1 & 1 & 1\\
\sin\alpha & 0 & -\sin\alpha & 0\\
0 & \sin\alpha & 0 & -\sin\beta_k\\
\cos\alpha & \cos\alpha & \cos\alpha & \cos\beta_k
\end{array}\right)
\end{equation}

\noindent and the superindex $\bot$ denotes matrix transposition.
Then the infimum in (\ref{eq9.43}) is bounded from below by the
smallest eigenvalue of $M M^\bot$. The characteristic equation of
$MM^\bot$ is
\begin{eqnarray}\nonumber
&&\hskip -1cm -\lambda^4+8\lambda^3-\frac 14
\Bigl(55-22\cos(2\alpha)-3\cos(4\alpha)-8\cos(\alpha-\beta_k)
-\cos(2\alpha-2\beta_k)-2\cos(2\beta_k)\\[4pt] \nonumber
&&\hskip -1cm  -16\cos(\alpha+\beta_k)-3
\cos(2\alpha+2\beta_k)\Bigr)\lambda^2- \frac 12
\sin^2\alpha\Bigl(-4\cos(2\alpha)+\cos(4\alpha)+12
\cos(\alpha-\beta_k)\\[4pt]
\nonumber &&\hskip -1cm
-33+\cos(2\alpha-2\beta_k)+20\cos(\alpha+\beta_k)
+3\cos(2\alpha+2\beta_k)\Bigr)\lambda
=4\sin^2\alpha(\cos\alpha-\cos\beta_k)^2.
\end{eqnarray}

\noindent By the Mean Value Theorem for the function $\arccos$ and
our assumptions on $\alpha,\beta_k$ there exists $C_0(\alpha)$
independent of $\beta_k$ such that for all $k$
\begin{equation}\label{eq9.45}
|\alpha-\beta_k|\leq C_0(\alpha)|\cos\alpha-\cos\beta_k|,
\end{equation}

\noindent and therefore,
\begin{equation}\label{eq9.46}
4\sin^2\alpha(\cos\alpha-\cos\beta_k)^2\geq 4\sin^2\alpha
(C_0(\alpha))^{-2} |\alpha-\beta_k|^2.
\end{equation}

\noindent It follows that
\begin{equation}\label{eq9.47}
\lambda\geq \frac{\sin^2\alpha (C_0(\alpha))^{-2}}{100}\,
|\alpha-\beta_k|^2,
\end{equation}

\noindent because otherwise the left-hand side of (\ref{eq9.45})
is strictly less than its right-hand side. Combined with
(\ref{eq9.43}), this finishes the proof of (\ref{eq9.35}). The
statement (\ref{eq9.36}) follows from (\ref{eq9.35}) and
Theorem~\ref{t1.2}.\ep

\noindent {\it Remark}. Retain the conditions of Lemma~\ref{l9.4}, and let $b:=a^{1/5}$.
By our construction, in each spherical
layer ${\overline{C_{b^{-j},b^{-j+1}}}}$ there are  either\\
(i) exactly three points $A_i^k$, $i=1,2,3$ for some
$k=k(j)$,\\
(ii) or exactly one point $A_4^k$, $k=k(j)$,\\
(iii) or no points from the complement of $\Omega$.

By Lemma~\ref{l9.2} it follows that in either case
\begin{equation}\label{eq9.48}
{{\rm
Cap}\,(\overline{C_{b^{-j},b^{-j+1}}}\setminus\Omega,\RR^3\setminus
\{O\})}=0
\end{equation}

\noindent and hence,
\begin{equation}\label{eq9.49}
\sum_{j} b^{-j} \,\inf_{P\in\Pi_1}{{\rm
Cap}_P\,(\overline{C_{b^{-j},b^{-j+1}}}\setminus\Omega)} =0.
\end{equation}

\noindent At the same time, if $\sum_{k}(\alpha-\beta_k)^2$
diverges, then so does the integral in (\ref{eq9.35}).

It follows that for the {\it same} domain $\Omega$ the convergence
of the integral in (\ref{eq1.9}) might depend on the choice of
$a$.

Alternatively, one can say that for the same $a$ the convergence
of the integral in (\ref{eq1.9}) might depend on the dilation of
the domain. In particular, (\ref{eq1.9}) {\it can not} be a
necessary condition for the 1-regularity since the concept of
1-regularity is dilation invariant.

Conversely, our proof of the first statement in Theorem~\ref{t1.2}
and Proposition~\ref{p7.1} relies on Proposition~\ref{l5.4} which,
in turn, follows from the Poincar{\'e}-type inequality
(\ref{eq5.21}). In fact, for every $s$ our choice of $P$, that
allows to estimate the infimum under the integral sign in
(\ref{eq1.9}), is dictated by the approximating constants in the
Poincar{\'e}'s inequality on $(s,as)$ (see the proof of
Lemma~\ref{l5.3}). Therefore, in our argument one can not make a
uniform choice of $P$ for all $s$ to substitute (\ref{eq1.9}) with
(\ref{eq1.10}).

\begin{corollary}\label{c9.5} The 1-irregularity may be unstable under the affine transformation of
coordinates.
\end{corollary}

\bp The proof is based on Corollary~\ref{c9.3} and
Lemma~\ref{l9.4}. Indeed, given the assumptions of
Lemma~\ref{l9.4}, if $\beta_k=\alpha$ for all $k$, then the
complement of $\Omega$ is entirely contained in the circular cone
of aperture $\alpha$ with the vertex at the origin and hence, by
virtue of Corollary~\ref{c9.3}, the point $O$ is not 1-regular.

However, if $\beta_k=\alpha+\eps$ for all $k$, then the series in
(\ref{eq9.36}) diverges for arbitrary small $\eps>0$, which
entails 1-regularity of $O$.\ep

\vskip 0.08in \noindent --------------------------------------
\vskip 0.10in

\noindent {\it Svitlana Mayboroda}

\noindent Department of Mathematics, Brown University,\\
151 Thayer Street, Providence, RI 02912

\vskip 0.05in

\noindent Department of Mathematics, The Ohio State University,\\
231 W 18th Av., Columbus, OH, 43210, USA

\noindent {\tt svitlana\@@math.ohio-state.edu,
svitlana@math.brown.edu}

\vskip 0.10in

\noindent {\it Vladimir Maz'ya}

\noindent Department of Mathematics, The Ohio State University,\\
231 W 18th Av., Columbus, OH, 43210, USA

\vskip 0.05in

\noindent Department of Mathematical Sciences, M\&O Building,\\
University of Liverpool, Liverpool L69 3BX, UK

\vskip 0.05in

\noindent Department of Mathematics, Link\"oping University,\\
SE-581 83 Link\"oping, Sweden

\noindent {\tt vlmaz\@@math.ohio-state.edu,  vlmaz\@@mai.liu.se}

\end{document}